\documentclass[preprint,12pt,3p]{elsarticle}

\usepackage{amssymb}
\usepackage{amsmath}
\usepackage{graphicx}
\usepackage{float}
\usepackage{caption}
\usepackage{subcaption}
\usepackage{import}
\usepackage{color}
\usepackage{makecell}
\usepackage{tikz}
\usepackage{lettrine}
\usepackage{placeins}
\usepackage{hyperref}

\definecolor{blau}{RGB}{0,0,0}
\definecolor{vermell}{RGB}{0,0,0}

\graphicspath{{./images/}}

\journal{CSMA 2017}

\begin{document}

\begin{frontmatter}

\title{Explicit topology optimization through moving node approach: beam elements recognition}

\author[label1]{Ghislain Raze\footnote{ Research Assistant, Structural Mechanics Department.}\corref{cor1}}
\address[label1]{Institut Cl\'ement Ader, Universit\'e de Toulouse, ISAE-SUPAERO-CNRS-INSA-Mines Albi-UPS, Toulouse, France}
\ead{ghislain.raize@isae-supaero.fr}

\author[label1]{Joseph Morlier\footnote{Professor, Structural Mechanics Department.}}
\ead{joseph.morlier@isae-supaero.fr}
\cortext[cor1]{Corresponding author}

\begin{abstract}
 Structural optimization (topology, shapes, sizing) is an important tool for facilitating the emergence of new concepts in structural design. Normally, topology optimization is carried out at the early stage of design and then shape and sizing design are performed sequentially. Unlike traditional topology optimization method, explicit methodologies have attracted a great deal of attention because of the advantages of shortcuting the costly CAD/CAE processes while dealing with low order number of design variables compared to implicit method (such as SIMP). This paper aims at presenting an adaptation of a flow-inspired approach so-called Moving Node Approach (MNA) in topology optimization. In this approach, the discretization is decoupled from the material distribution and the final objective is to recognize the best beam assembly while minimizing compliance. The paradigm has here changed and new design variables are used such as nodes location, elements length/orientation and size providing a lower number of design variables than pixels-based. The methodology is validated using 2 classical testcases in the field of Topology Optimization: the Cantilever beam and the L-Shape problem.
\end{abstract}


\end{frontmatter}

\section{Introduction}
Advanced shape and topology optimization methods have been addressed as the most promising techniques for least-weight and performance design of engineering structures. Since the pioneering work of \cite{bendsoe1988generating}, topology optimization  has received considerable research attention. Numerous topology optimization approaches such as SIMP (Solid Isotropic Material with Penalization) approach \cite{bendsoe1989optimal,zhou1991coc}; , level set approach \cite{wang2003level,allaire2004structural} and evolutionary approach \cite{xie1993simple} have been studied. We can simply define topology optimization as the process to find the best material layout in a domain in order to maximize specific performance targets. In the so-called implicit topology optimization methods, each element is associated with a density indicating the presence or absence of material. These densities can either be discrete or continuous and serve as optimization variables. The level-set method (\cite{sethian2000structural}) is an interesting approach in which the material layout is described as a level-set function.
It has been extent to large problem \cite{suresh2013efficient} or for educational purposes \cite{suresh2010199} using Pareto methodology. Some problems arose from these techniques. First, the presence of intermediate densities can be troublesome for the part's manufacturing. It can be penalized by considering that the material's stiffness does not scale linearly with the density, but with a power law. The power which appears to work best and which is commonly used is $p=3$. Second, the numerical solution can follow a checkerboard pattern. A filter inspired from image processing was proposed by \cite{sigmund1994design} in order to avoid holes under a given length scale. Normally, topology optimization is carried out at the early stage of design process. Then a new costly design cycle must be relaunched: mesh creation from topology optimization results, and if validated, a new sizing optimization problem.  Some researchers try then to extract directly from topology optimization results (image processing) the structure skeleton to relaunch easily the FE sizing process \cite{yi2017identifying,gedig2010framework}. Moreover as demonstrated by \cite{guo2014topology}, the number of design variables involved in implicit topology optimization approaches is relatively large especially for three dimensional problems. Recently, Liu proposed a way to reduce the dimensionality of implicit topology using machine learning clustering method \cite{liu2015thin}. The same year, Gogu introduced the use reduced order models for improving the efficienty of large problem \cite{gogu2015improving}. 
To authors's knowledge the very first work in explicit topology optimization is related to geometrical projection's method \cite{bell2012geometry} derived from shape optimization author's previous works \cite{norato2004geometry}.Recently a so-called moving morphable components based topology optimization framework, has been developed in \cite{guo2014topology}. The originality of this method is to build blocks of topology optimization through a set of morphable components. Recently, the same idea has also been adopted for beams and plates (\cite{norato2015geometry,zhang2016geometry}) based on the SIMP framework for topology optimization of continuum structures made of discrete elements.

This work explores a flow-inspired method so-called Moving Node Approach (MNA)\cite{overvelde2012moving} to find the optimal structural topologies by optimizing the shapes, lengths, thicknesses, orientations and layout (connectivity) of these components. Methodology of adaptation/derivation of MNA are investigated in detail in section 2. By employing the position of the nodes as design variables in the topology optimization method. The topology optimization problem then transforms into a flow-like problem, in which the material moves to a more optimal distribution. We first compare classical linear elasticity Finite Element analysis (FE) to the so-called Element-Free Galerkin (EFG) method which was proposed by Belytschko \cite{belytschko1994element}. We compare the two approaches on the well-know cantilever beam problem.  In addition to MNA, a modification is done to both reduce the complexity and obtain an explicit final design. We develop merging procedure during the iteration in order to do in the loop structural element recognition. We validate this procedure on the L-shape Case. The procedure for the MNA-based algorithm is simple, practical, requires little expert interference.  The performance of the algorithm is validated by these two classical examples. The authors want to emphasis that this method can be used as a conceptual design tool, that's why we propose some good practice and use often coarse mesh. The results show it can accurately extract structural element on the 2 testcases. Our FE code (topmna.m) is inspired from \cite{andreassen2011efficient}, itself inspired by \cite{sigmund200199} and recently extend to 3D \cite{liu2014efficient}.

\section{The Moving Nodes Approach framework}

    \label{sec:methodology}
    \subsection{Motivation and basic idea}
    
        Unlike traditional Topology Optimization approaches, explicit approaches provide a geometrical description of the structure, using primary building blocks, often referred to as components. These components can move change shape and dimensions, and also overlap with each other. The structural topology is represented by the union of the components. Explicit approaches by nature have the merit of providing black and white designs, this is generally due to the formulation of density field of a component, which does not allow intermediate densities, see figure 1 for reference:
        
        \begin{figure}[ht]
            \centering
            \vspace{-0.1cm}
            \includegraphics[width=0.7\textwidth]{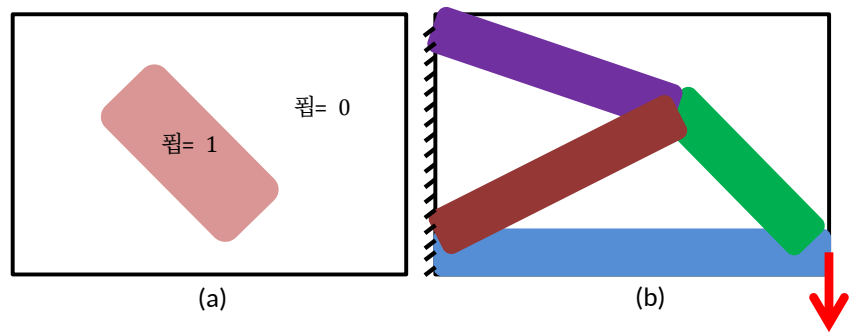}
            \vspace{-0.1cm}
            \caption{(a)-density field of a component, (b) structure's topology using components}
            \label{figure1}
        \end{figure}
        
        However, it appears in practice that explicit approaches involving components with such a topology and density field present some difficulties to converge, as will be explained in detail latter. In addition, the final design depends strongly on the initial design and in some cases the optimization may not even converge for a disconnected initial configuration. In order to tackle these challenges we propose a new geometric approach for explicit topology optimization, called the Moving Nodes Approach (MNA).
        
    \subsection{Geometry description and density field}
    
        The moving nodes approach is originally inspired from the work of \cite{overvelde2012moving}, where a meshless and flow-inspired method for topology optimization was developed, using the Element-Free Galerkin method (EFG) instead of finite elements method (FEM). The approach we propose here consist on representing the structure with a (finite) number of components which centers are traditionally called mass nodes, each mass node has an influence region, it is the region occupied by the component. Density is equal to one in the mass nodes location, and decreases to zero at the borders of the component. Here rectangular components are used and each component is described by following geometric design variables (same as in MMC method):
        
        $\left(x_0,y_0\right)\ $coordinates of the component's center (mass node)

        $L\ $half length of the component (along x axis)

        $t\ $half thickness of the component (along y axis)

        $\theta \ $orientation angle with respect to x axis

        Left figure shows a component's geometry and its design variables, and right one shows the corresponding density field:
    
        \begin{figure}[ht]
            \centering
            \vspace{-0.1cm}
            \includegraphics[width=0.7\textwidth]{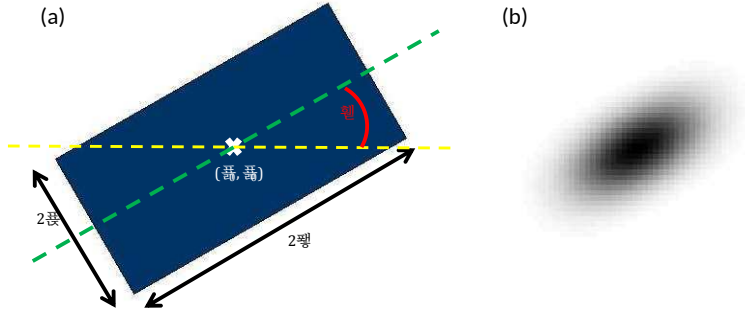}
            \vspace{-0.1cm}
            \caption{(a)-design variables of MNA component, (b)-density field of MNA component}
            \label{figure2}
        \end{figure}
        
        \begin{equation}
            \rho(\mathbf{x}) = \sum_{I=1}^n m^I W(\mathbf{x},\boldsymbol{\mu}^I) \label{eq:density}
        \end{equation}
        
        where $n$ denotes the number of mass nodes whose domain of influence covers point $\mathbf{x}$, $m^I$ is the mass associated to mass node $I$ and $\boldsymbol{\mu}^I$ is the associated material variables vector. In two dimensions and for a rectangular domain of influence, each mass node can be characterized by a position $(x^I,y^I)$, an orientation $\theta^I$ and domain dimensions $L_x^I$ and $L_y^I$. Similarly to the work of \cite{guo2014topology}, these variables can be modified to diminish the compliance.
        
        \begin{equation}
            \boldsymbol{\mu}^I = [x^I, y^I, \theta^I, L_x^I, L_y^I]^T
        \end{equation}
    
        If we consider the local variables
        
        \begin{equation}
            \begin{array}{c}
            \displaystyle \xi(\mathbf{x},\boldsymbol{\mu}^I) = \frac{(x-x^I)\cos(\theta^I) + (y-y^I)\sin(\theta^I)}{L_x^I/2} \\
            \displaystyle \eta(\mathbf{x},\boldsymbol{\mu}^I) = \frac{-(x-x^I)\sin(\theta^I) + (y-y^I)\cos(\theta^I)}{L_y^I/2}
            \end{array}
            \label{eq:localVar}
        \end{equation}
        
        the kernel function in two dimensions can be expressed as
        
        \begin{equation}
            W(\mathbf{x},\boldsymbol{\mu}^I) = w(|\xi(\mathbf{x},\boldsymbol{\mu}^I)|,\frac{L_x^I}{2})w(|\eta(\mathbf{x},\boldsymbol{\mu}^I)|,\frac{L_y^I}{2})
            \label{eq:shape}
        \end{equation}
        
        where $|x|$ is the modulus of $x$ and $w$ is the kernel function in one dimension. It can be chosen as the cubic spline weight function
        
        \begin{equation}
        \def\arraystretch{2}
        w(r,d) = \left\{ \begin{array}{l l}
             \displaystyle\frac{2}{d}(\frac{2}{3} - 4 r^2 + 4 r^3) \quad & r \leq \frac{1}{2}  \\
             \displaystyle\frac{2}{d}(\frac{4}{3} - 4r + 4r^2 - \frac{4}{3} r^3) \quad & \frac{1}{2} < r \leq 1 \\
             \displaystyle 0 & r > 1
        \end{array}\right.
        \label{eq:cswf}
    \end{equation}
    
    The constant $\frac{2}{d}$ is set so that the integral of the kernel function over its whole domain of influence is equal to one. The number $d$ corresponds to the smoothing length of the kernel function.
    
    \subsection{Density derivatives}
    
        In order to know how the material distribution affects the compliance, computation of the density derivatives is required. Let $x_i^I$ be any material variable of mass node $I$. The general expression of the derivative of the density with respect to that variable is
    
        \begin{equation}
            {\rho(\mathbf{x}^k)}{x_i^I} = {m^I}{x_i^I}W(\mathbf{x},\boldsymbol{\mu}^I) + m^I{W(\mathbf{x},\boldsymbol{\mu}^I)}{x_i^I} \label{eq:densityDerivatives}
        \end{equation}
        
        The mass of mass node $I$ is proportional to its domain dimensions ($\beta$ is the density, constant here).
        
        \begin{equation}
            m^I = \beta L_x^I L_y^I \label{eq:nodalMass}
        \end{equation}
        
        Using \eqref{eq:localVar}, \eqref{eq:shape}, \eqref{eq:cswf} and \eqref{eq:nodalMass}, the general expression \eqref{eq:densityDerivatives} can then be computed analytically for any material variable. \\
        
        Now let three cases be defined:
        
        \begin{itemize}
            \item If only $x^I$ and $y^I$ are modified, $I$ refers to a \textit{mass node}.
            \item If $x^I$, $y^I$ and $\theta^I$ are modified, $I$ refers to mass node specialized into an \textit{undeformable structural member}.
            \item If all the material variables can be modified, $I$ refers to a mass node specialized into a \textit{deformable structural member}.
        \end{itemize}
    
    \subsection{Compliance sensitivity}
    
        The compliance is the work done by external forces and can be considered as the inverse of the global structural stiffness. Therefore, it is the objective function that should be minimized in the optimization algorithm. It is given by the scalar product of the nodal force vector $\mathbf{F}$ and the nodal displacements $\mathbf{U}$.
    
        \begin{equation}
            C = \mathbf{F}^T \mathbf{U} \label{eq:compliance}
        \end{equation}
    
        Taking the derivative of this expression with respect to $x_i^I$ yields
    
        \begin{equation}
            {C}{x_i^I} = {\mathbf{F}}{x_i^I}^T\mathbf{U} + \mathbf{F}^T{\mathbf{U}}{x_i^I} \label{eq:dcompliance}
        \end{equation}
    
        The considered problems are statically loaded and without body forces, therefore the nodal force vector does not depend on the density distribution and the first term in the right hand side of \eqref{eq:dcompliance} is equal to zero. Let us now consider the discrete equilibrium equation \eqref{eq:equilibrium}.
    
        \begin{equation}
            \mathbf{K}\mathbf{U} = \mathbf{F}\label{eq:equilibrium}
        \end{equation}
        
        where $\mathbf{K}$ is the stiffness matrix. The derivative of this equation with respect to $x_i^I$ is
        
        \begin{equation}
            {\mathbf{K}}{x_i^I}\mathbf{U} + \mathbf{K}{\mathbf{U}}{x_i^I} = {\mathbf{F}}{x_i^I} = 0 \label{eq:dudxi}
        \end{equation}
        
        Inserting \eqref{eq:equilibrium} and \eqref{eq:dudxi} into \eqref{eq:dcompliance} and remembering that the stiffness matrix is symmetric yields the final expression of the compliance sensitivity
        
        \begin{equation}
            {C}{x_i^I} = -\mathbf{U}^T {\mathbf{K}}{x_i^I} \mathbf{U}
        \end{equation}
        
        Therefore, only the derivatives of the stiffness matrix with respect to the material distribution variables are to be evaluated. The stiffness matrix is assembled with a Gauss quadrature
        
        \begin{equation}
            \mathbf{K} = \sum_{k=1}^{n_G} \mathbf{K}_e^k E(\mathbf{x}^k) \omega^k 
        \end{equation}
        
        where $\mathbf{x}^k$ is the coordinates vector of Gauss point $k$, $\omega^k$ is its associated weight, $E$ is the Young modulus and $\mathbf{K}_e^k$ is its associated element stiffness matrix with a unit Young modulus. It is computed thanks to a discretization technique: with the FEM, it is integrated with the Gauss points of the elements, while with the EFG method it is computed with the Gauss points of the background mesh integration cells. The Young modulus $E$ is the only quantity depending on the material distribution variables.
        
        \begin{equation}
            {\mathbf{K}}{x_i^I} = \sum_{k=1}^{n_G} \mathbf{K}_e^k {E(\mathbf{x}^k)}{x_i^I} \omega^k 
        \end{equation}
        
        The Young modulus itself directly depends on the density function. This will be explained hereafter with equations \eqref{eq:PLERhoMin} and \eqref{eq:PLEMin}. 
    
    \subsection{Optimizer details}
    
            
            Gradient-based optimizers are usually fast but can converge to local optima. In the scope of this work, a steepest descent algorithm, a conjugated gradients algorithm and a quasi-Newton BFGS algorithm were implemented (details about these algorithms can be found in \cite{nocedal2006numerical} or \cite{craveur2014optimisation}). The Matlab's functions \texttt{fminunc()} and \texttt{fmincon()} also allow to proceed to unconstrained and constrained gradient-based optimization.
            

    \subsubsection{Mass constraint}
            
            If deformable structural members are used, they have the tendency to grow larger in size and make bulky structure. To avoid this and remain with a fixed maximum mass, a constraint on the total mass can be used, based on the material variables
            
            \begin{equation}
                \Big(m_{Max} - \sum_{I=1}^N \beta L_x^I L_y^I\Big) \geq 0
                \label{eq:massConstraint}
            \end{equation}
    
           $m_{Max}$ is the maximum allowed structural mass. The mass effectively used for structural stiffness is less than or equal to the mass given by summing all the nodal masses, since the influence domain of the mass nodes can be partly uncovered by the mesh. The structural mass can be obtained by integrating the density field over the whole mesh. The evaluation of the constraint and its gradient would however be less direct. Using the constraint \eqref{eq:massConstraint} allows to project easily any layout on the admissible domain since it is a quadratic function of the material variables.
    
    \subsection{Numerical aspects}
        
            The mass nodes have a natural tendency to stack on top of each other, resulting in zones where the density is greater than one. From a manufacturer's point of view, crafting the part then becomes difficult, if not impossible. Therefore a density which ranges from zero to one is required. As in \cite{overvelde2012moving}, one can use the asymptotic density
            
            \begin{equation}
                \rho^a(\mathbf{x}) = \frac{a\rho(\mathbf{x})}{(\rho(\mathbf{x}))^b+a}
            \end{equation}
            
            and change the density derivatives accordingly.
            
            \begin{equation}
                \frac{d\rho^a(\mathbf{x})}{d\rho(\mathbf{x})} = \frac{a(1-b)(\rho(\mathbf{x}))^b + a^2}{[(\rho(\mathbf{x}))^b+a]^2}
            \end{equation}
            
            with
            
            \begin{equation}
                b = \frac{1}{1-\rho_{Max}}-1 \qquad \qquad a = \frac{(\rho_{Max})^b}{\rho_{Max}-1}
            \end{equation}
            
            The asymptotic density is almost linear when $0\leq\rho(\mathbf{x}) \leq \rho_{Max}$, and densities above $\rho_{Max}$ are strongly penalized. Hence $\rho_{Max}$ should be slightly greater than one.
            
            \begin{figure}[ht]
                \centering
                \includegraphics[width=0.325\textwidth]{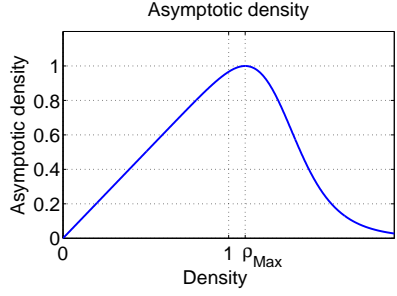}
                \vspace{-0.3cm}
                \caption{Example of asymptotic density with $\rho_{Max} = 1.2$}
                \label{fig:Asymptotic density}
            \end{figure}
        
            The asymptotic density thus naturally avoids designs with densities greater than one. It can however induce strong variations on short length scales, especially if $\rho_{Max}$ is close to one.
        
        \subsubsection{Filtering}
            
            In the SIMP, a filter can be used to avoid checkerboard patterns. A similar concept can be defined in the MNA by designing a filter based on the Gauss points. From Gauss points at coordinates $\mathbf{x}_i$ and $\mathbf{x}_j$, the convolution matrix can be defined
            
            \begin{equation}
                H_{kl} = H(\mathbf{x}_k,\mathbf{x}_l) = r_{Min} - \min(r_{Min},||\mathbf{x}_k-\mathbf{x}_l||)
            \end{equation}
            
            and modify the densities $\rho$ to filtered densities $\hat{\rho}$
        
            \begin{equation}
                \displaystyle\hat{\rho}(\mathbf{x}_{k}) = \frac{1}{\displaystyle \sum_{l=1}^{n_G} H_{kl}}\sum_{k=1}^{n_G} H_{kl}\rho(\mathbf{x}_l)
            \end{equation}
            
            The densities derivatives are filtered similarly.
        
        \subsubsection{Minimum density/minimum stiffness}
            
            Regions with zero density will have zero associated stiffness. This will result in a singular stiffness matrix. To avoid this, a minimum density $\rho_{Min}$ can be used.
            
            \begin{equation}
                \rho(\mathbf{x}) = \rho_{Min} + (1-\rho_{Min})\sum_{I=1}^n m^I W(\mathbf{x},\boldsymbol{\mu}^I)
                \label{eq:rhoMin}
            \end{equation}
            
            Since the density undergoes several transformations, a minimum Young modulus $E_{Min}$ can be used instead.
            
            \begin{equation}
                E(\mathbf{x}) = E_{Min} + (E_{0}-E_{Min})\rho(\mathbf{x})
                \label{eq:EMin}
            \end{equation}
            
            where $E_0$ is the material's Young modulus.

        \subsubsection{Intermediate density penalization}
        
            The SIMP intermediate density penalization can also be used. If a minimum density is used as in \eqref{eq:rhoMin}, then the Young modulus is given by
        
            \begin{equation}
                E(\mathbf{x}) = E_{0} (\rho(\mathbf{x}))^p \label{eq:PLERhoMin}
            \end{equation}
            
            Otherwise if a minimum Young modulus is used as in \eqref{eq:EMin}, then
            
            \begin{equation}
                E(\mathbf{x}) = E_{Min} + (E_{0}-E_{Min})(\rho(\mathbf{x}))^p \label{eq:PLEMin}
            \end{equation}
        
            The goal of this density penalization in the SIMP approach is to avoid as much as possible elements with intermediate density in the final design. Of course, this will not suppress intermediate densities in the MNA, but this gives a common ground to compare both methods.
            
     \subsubsection{Proposed algorithm and element merging}
     
\begin{figure}[ht]
                \centering
                \includegraphics[width=0.75\textwidth]{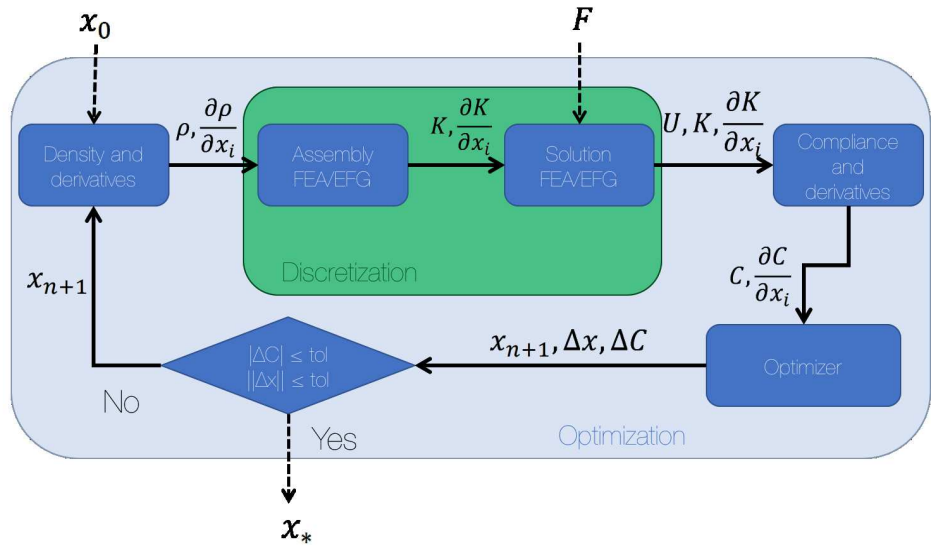}
                \vspace{-0.3cm}
                \caption{Proposed algorithm: the user can use either FEA or EFG, adjust tolerence to converge to a final design. $x_0$ is here a vector of size N*5 where N is the number of beam element. If merging procedure is used $x*$ can be reduced to the sufficient number of elements (beam assembly) needed to minimize the compliance}
                \label{fig:bp}
            \end{figure}

            To resume our methodology we propose a scheme which is depicted in figure \ref{fig:bp}. MNA offers the possibility to use Mass Nodes (2 variables per element: $x$ and $y$), undeformable element (3 variables per element: $x$, $y$, $\theta$ and deformable element (5 variables per element:$x$, $y$, $\theta$, $L_x$, $L_y$). So the initial size of the deformable element problem $x_0$ is normally regularly sampled on the domain and has a size equal to $N*5$ with N the number of deformable element (beam). On can notice that if merging procedure is used $x*$ can be reduced to the sufficient number of elements (beam assembly) needed to minimize the compliance.

                 Although the material distribution and the discretization are decoupled, they should not be chosen independently. Moreover, the aforementioned possibility of strong variations due to the asymptotic density can induce a wrong integration. The following rule of thumb allowed the authors to avoid any problem with a density $\rho_{Min} = 1.05$
        
        \begin{equation}
            \min_I (d_x^I,d_y^I) \geq \max_K (\Delta x^K, \Delta y^K) 
            \label{eq:goodPractice}
        \end{equation}
        
        Where $\Delta x^K$ is the size along $x$ of element/integration cell $K$. This means that a mass node influence domain should cover at least two elements/integration cells in one direction. When deformable structural components are used, this rule does not suffice in general, as $d_x^I$ and $d_y^I$ can change. This issue will be discussed later.
            
\begin{figure}[ht]
                \centering
                \includegraphics[width=0.4\textwidth,clip=true,trim=0cm 1cm 0cm 1cm]{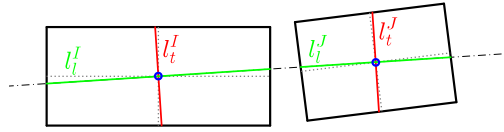}
                \caption{Merging procedure}
                \label{fig:merging}
            \end{figure}

            After convergence, similar nodes could be merged. First, their orientations are compared. If they are equal up to a certain tolerance, their dimensions are compared and the distance between their two centers is compared to their dimensions. If the criteria are met, the nodes are merged. Consider the two nodes shown in figure \ref{fig:merging}. They are merged if they satisfy the following conditions :
             
             \begin{enumerate}
                 \item $\displaystyle|\theta^I-\theta^J| \leq tol_\theta$, where $\displaystyle tol_\theta$ is a tolerance on orientations
                 \item $\displaystyle|l_t^I - l_t^J| \leq tol_l$, where $\displaystyle tol_l$ is a tolerance on lengths
                 \item $\displaystyle || \mathbf{x}^I - \mathbf{x}^J|| \leq (1+tol_d) r_\rho \frac{1}{2}(l_l^I + l_l^J)$, where $\displaystyle tol_d$ is a tolerance on distances and $\displaystyle r_\rho $ is a ratio related to the density level
             \end{enumerate}
             
            The merging step in the MNA can be related to the filtering step in the SIMP, although it has some differences. It does not suppress small mass nodes, and it is not done after each iteration but rather after convergence (and the optimization algorithm starts again with the simplified structure).


\section{Numerical experiments}\label{Section:Experiment}

\subsection{FE, EEG comparison on Cantilever beam testcase}
\begin{figure}[ht]
            \centering
            \includegraphics[width=0.25\textwidth]{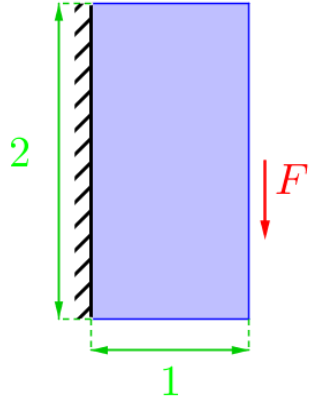}
            \caption{Test case}
            \label{fig:testCase}
        \end{figure}

        The test case shown in figure \ref{fig:testCase} is analyzed. The designable space dimensions are $1\times2$. The frame of reference is set so that $(x,y)\in[0,1]\times[-1,1]$. It is clamped at its left boundary $x = 0$ and loaded by a unit force at $(x,y) = (1,0)$.
        
        The other important quantities are listed in table \ref{tab:testCase}.
        
        \begin{table}[ht]
            \footnotesize
            \def\arraystretch{1.2}
            \centering
            \begin{tabular}{|c|c|c|}
                \hline
                \textbf{Quantity} & \textbf{Symbol} & \textbf{Value} \\
                \hline
                 Young modulus & $E_0$ & 1  \\
                \hline
                 Poisson ratio & $\nu$ & 0.3 \\ 
                \hline
                Elements/cells per unit length &$n_x$ or $n_y$ & [4 ; 15] \\
                \hline
                Gauss points per element/cell & $n_{G,e}$ & 4 \\
                \hline
                Shape functions degree & $p$ & 1 \\
                \hline
                Num. of mass nodes & $n_{mn}$ & 4 \\
                \hline
                Volume fraction & $v_{Frac}$ & [0.05 ; 1] \\
                \hline
                Max. iteration & $iter_{Max}$ & $1,000$ \\
                \hline
                Max. tol. on compliance change & $tol_C$ &$1e^{-6}$\\
                \hline
                Max. tol. on variables change & $tol_x$ & $1e^{-6}$\\
                \hline
                Max. tol. on mass constraint & $tol_m$ & $1e^{-6}$\\
                \hline
                Discr. relative smoothing length & $d$ & 2.5\\
                \hline 
                Mass relative smoothing length & $d^\rho$ & 1.5 \\
                \hline
            \end{tabular}
            \caption{Test case numerical values}
            \label{tab:testCase}
        \end{table}
        
        To avoid as much as possible making assumptions about the final shape, deformable structural members are used unless stated otherwise.

        This paragraph aims at showing the effect of discretization/optimizers on the final results, the computational complexity/memory of all tested methods,and finally volume fraction and filtering influence.
        
        The 88 lines program by \cite{andreassen2011efficient} is used to compare the results.
        
        The effects of the discretization method, the optimizer, the volume fraction and the filter are studied.
    \subsection{Optimizer effect}
        
\begin{figure}[ht]
            \centering
            \vspace{-0.4cm}
            \includegraphics[width=0.35\textwidth]{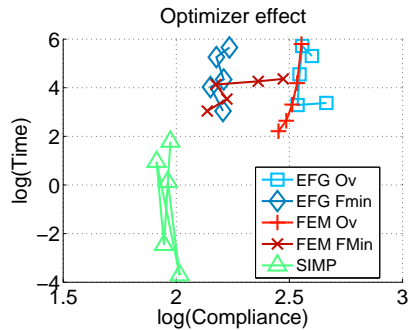}
            \vspace{-0.25cm}
            \caption{Optimizer effect for the EFG-based MNA and the FEM-based MNA: Overvelde's algorithm (Ov) and Matlab's \texttt{fmincon()} (FMin)}
            \label{fig:optimizer}
        \end{figure}

        Compliances and CPU times are compared for the MNA and the SIMP. The maximum volume fraction is constant and set to $v_{Frac} = 0.33$ while the discretization parameters $n_x$ and $n_y$ vary in their intervals. Three different optimizers for the MNA have been tested: the one proposed by Overvelde \cite{overvelde2012moving} with a decreasing time step to limit the oscillations, Matlab's \texttt{fmincon()} (gradient-based optimizer with constraints) and Matlab's \texttt{ga()} (genetic algorithm). The latter required tremendous amounts of time and did not converge to satisfactory results and its results are therefore not displayed. The figure \ref{fig:optimizer} shows the results obtained for the other different approaches.
        
        First, the genetic algorithms are observed to be generally slower than the other algorithms, which is not surprising given the huge number of evaluations required. Their resulting compliance values are higher, which means that their final configuration is not as well optimized as the one given by the gradient-based algorithms.

        Second, the EFG-based MNA is much slower than the FEM-based MNA. This can be explained by several factors. The nodes in the EFG method usually have more neighbouring nodes than those in the FEM. Consequently, the assembly of the stiffness matrix and its derivatives is longer. Moreover, the bandwidth of these matrices is larger. Finally, the addition of Lagrangian multipliers destroys their structure, making the matrix inversion longer.
        
\begin{figure}[ht]
            \center
            \includegraphics[width=0.08\textwidth,clip=true,trim=5cm 1cm 4cm 1cm]{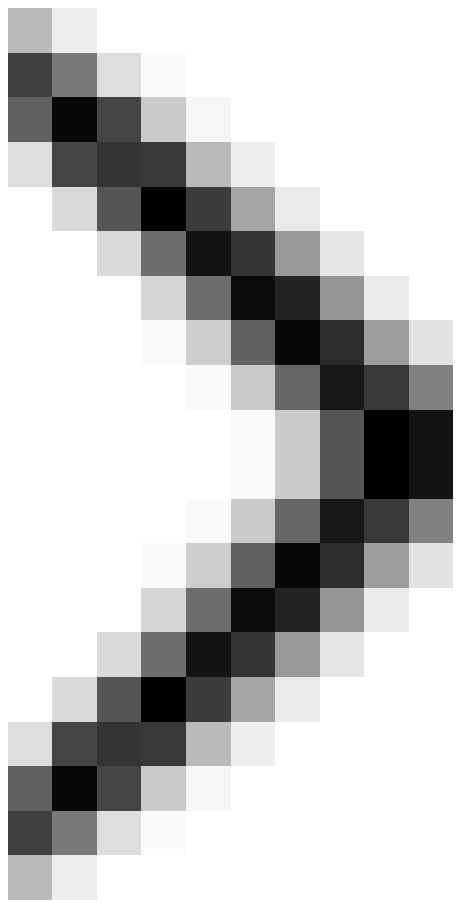}
            \includegraphics[width=0.08\textwidth,clip=true,trim=5cm 1cm 4cm 1cm]{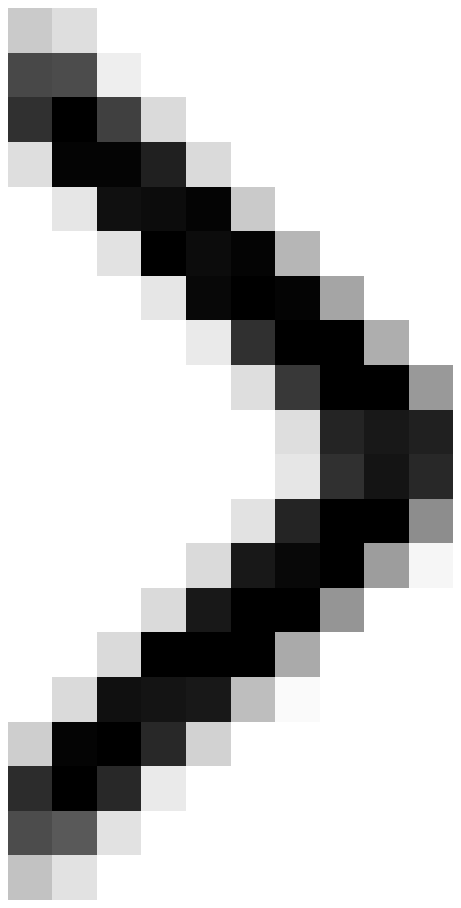}\\
            \includegraphics[width=0.08\textwidth,clip=true,trim=5cm 1cm 4cm 1cm]{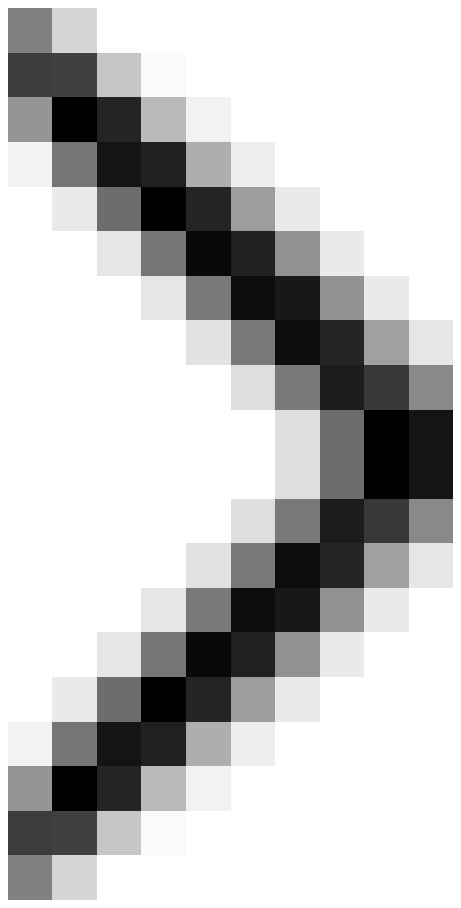}
            \includegraphics[width=0.08\textwidth,clip=true,trim=5cm 1cm 4cm 1cm]{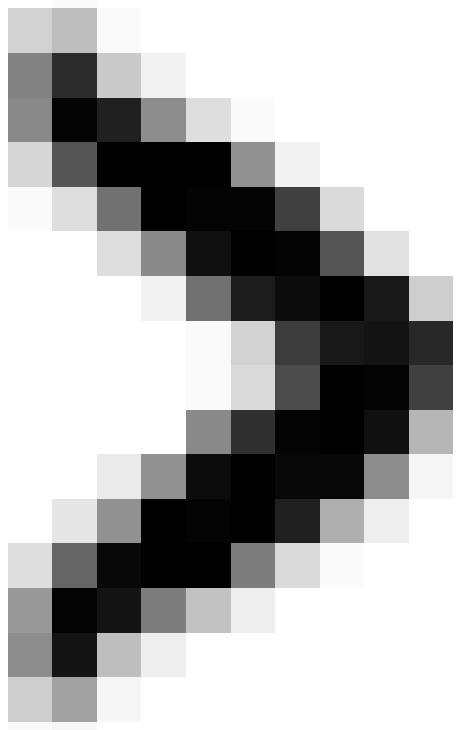}\\
            \vspace{-0.25cm}
            \caption{Density distribution for the EFG-based MNA (top) and the FEM-based MNA (bottom) for different optimization algorithms with 12 elements/cells per unit length: Overvelde's algorithm (left) and Matlab's \texttt{fmincon()} (right)}
            \label{fig:density1}
        \end{figure}

        The gradient-based optimizers therefore seem to be more adapted to the MNA even when the number of design variables is very low. The figure \ref{fig:density1} shows the final configurations obtained with the finest discretization.
        
    \subsection{Computational complexity}
        
        The computational complexity of the EFG-based MNA and the FEM-based MNA are compared with the SIMP in terms of computation time and memory used. The computations were done on a Intel Core i5-2410M CPU 2.30 GHz with 4 Gb of RAM.
        
\begin{figure}[ht]
            \centering
            \includegraphics[width=0.4\textwidth]{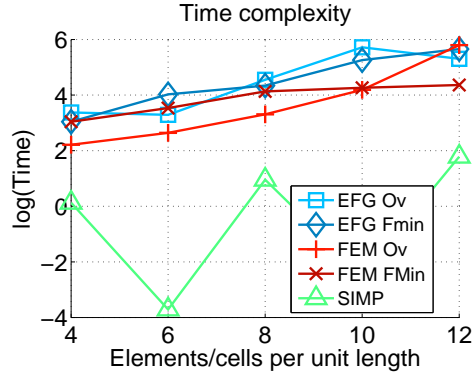}
            \vspace{-0.25cm}
            \caption{Time complexity}
            \label{fig:time}
        \end{figure}

        The figure \ref{fig:time} shows the time complexity of the different optimizers. Except for the SIMP, the MNA algorithms grow similarly in time as the discretization is refined.
        
        In terms of time, the cheapest discretization technique is the FEM. The most expensive algorithm is the genetic algorithm, followed by the Overvelde's algorithm and the Matlab's gradient-based algorithm. 
        
        The SIMP is always cheaper than the MNA. The execution time with 6 and 10 elements per unit length is quite larger than with 4, 8 or 12 elements. This is due to oscillations in the objective function which are not resolved before the maximum number of iterations is reached.
        
        \begin{figure}[ht]
            \centering
            \vspace{-0.4cm}
            \includegraphics[width=0.35\textwidth]{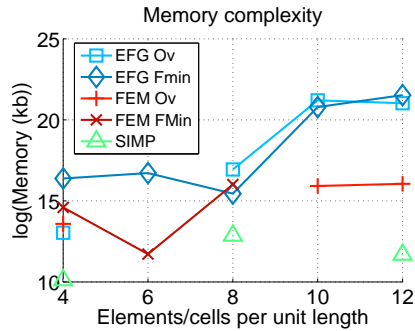}
            \vspace{-0.25cm}
            \caption{Memory complexity}
            \label{fig:memory}
        \end{figure}
        
        The memory complexity is displayed in the figure \ref{fig:memory}. The EFG-based MNA is once again the most expensive, due to the nature of its stiffness matrix and derivatives (non sparse). If there are no oscillations, the memory used by the SIMP is so small that it is set to zero by Matlab, hence not appearing on the logarithmic scale.
        
        The FEM-based MNA with an efficient gradient-based optimizer often seems to be the best in terms of time and memory complexity. It remains however a more time consuming than the SIMP. Indeed, the densities of the elements are the optimization variables in the SIMP. Hence, the assembly is much faster. Moreover, the optimality criteria method (described in \cite{sigmund200199}) is cheap and allows relatively fast convergence.
        
    \subsection{Volume fraction influence}
    
        The volume fraction influence is now studied for the EFG-based MNA, the FEM-based MNA and the SIMP. The volume fraction in the MNA cannot exactly be set, as a part of the total mass can flow out of the designable domain, and thus can not be taken into account in the structural stiffness. One can only set an approximate value of the volume fraction to begin with, and then get its evaluation by integrating the density over the domain as an output.
        
        \begin{figure}[ht]
            \centering
            \includegraphics[width=0.4\textwidth]{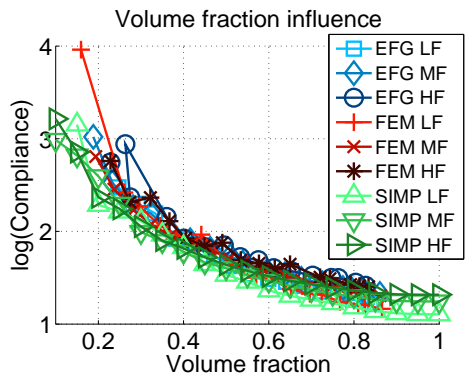}
            \vspace{-0.25cm}
            \caption{Volume fraction influence for the EFG-based MNA, the FEM-based MNA and the SIMP}
            \label{fig:vfrac}
        \end{figure}
        
        The tests are made with a Low Fidelity mesh (LF, $5\times10$ elements), a Medium Fidelity mesh (MF, $10\times20$ elements) and a High Fidelity mesh (HF, $15\times30$ elements). The optimizer used for the MNA is Matlab's \texttt{fmincon()}.
        
\begin{figure}[ht]
            \centering
            \includegraphics[width=0.08\textwidth,clip=true,trim=5cm 1cm 4cm 1cm]{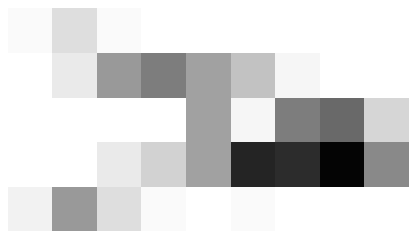}
            \includegraphics[width=0.08\textwidth,clip=true,trim=5cm 1cm 4cm 1cm]{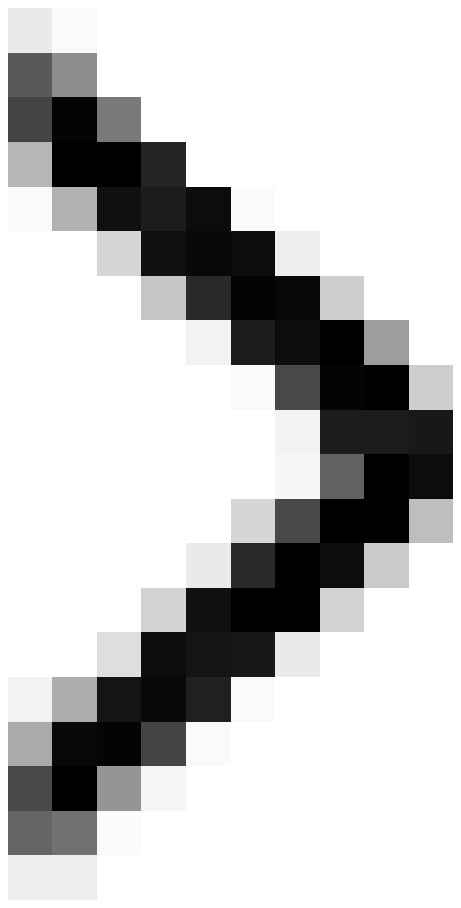}
            \includegraphics[width=0.08\textwidth,clip=true,trim=5cm 1cm 4cm 1cm]{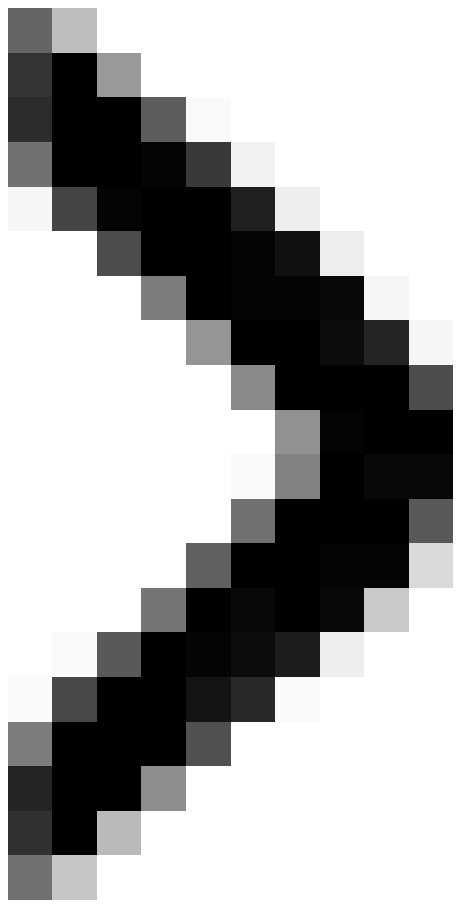}
            \includegraphics[width=0.08\textwidth,clip=true,trim=5cm 1cm 4cm 1cm]{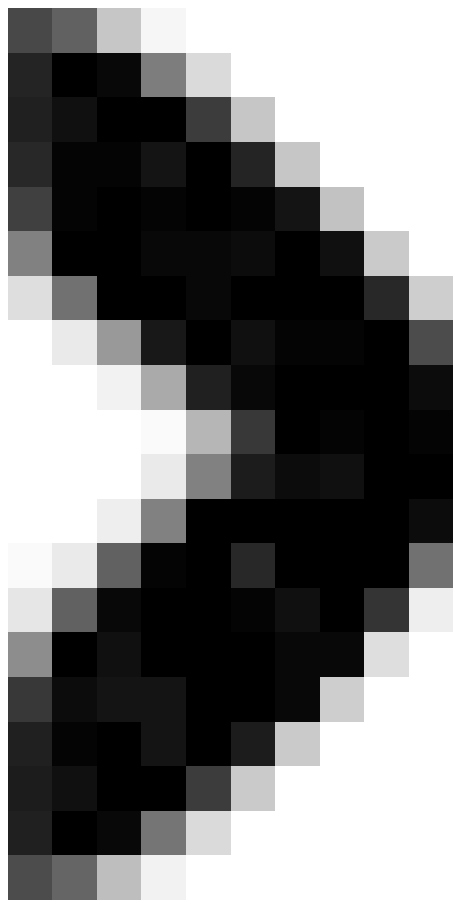}
            \includegraphics[width=0.08\textwidth,clip=true,trim=5cm 1cm 4cm 1cm]{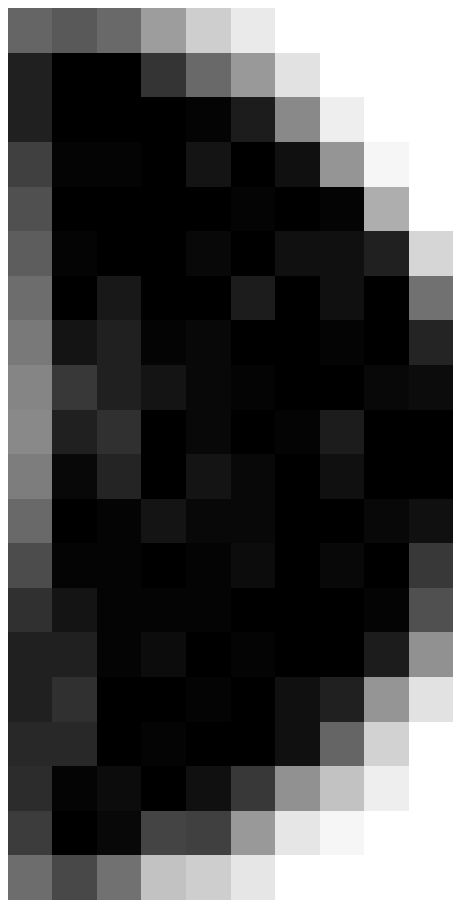}\\
            \includegraphics[width=0.08\textwidth,clip=true,trim=5cm 1cm 4cm 1cm]{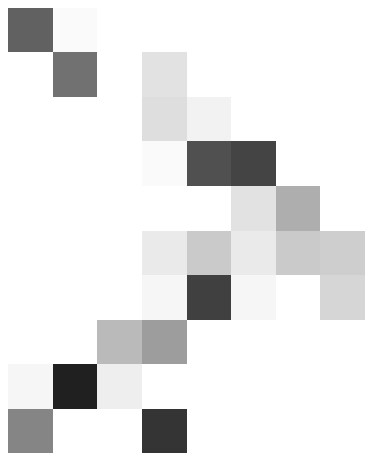}
            \includegraphics[width=0.08\textwidth,clip=true,trim=5cm 1cm 4cm 1cm]{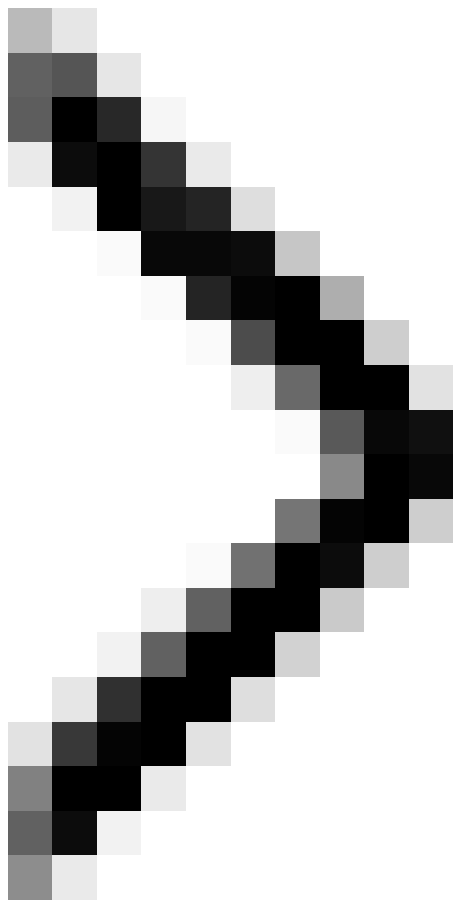}
            \includegraphics[width=0.08\textwidth,clip=true,trim=5cm 1cm 4cm 1cm]{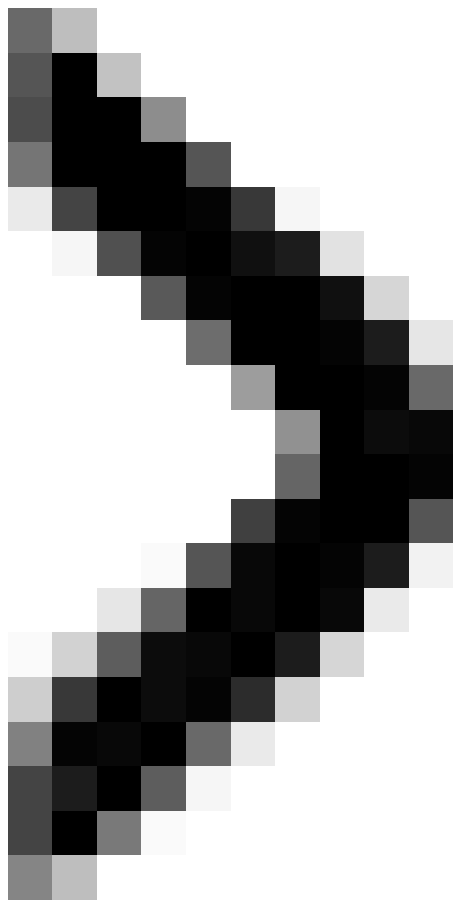}
            \includegraphics[width=0.08\textwidth,clip=true,trim=5cm 1cm 4cm 1cm]{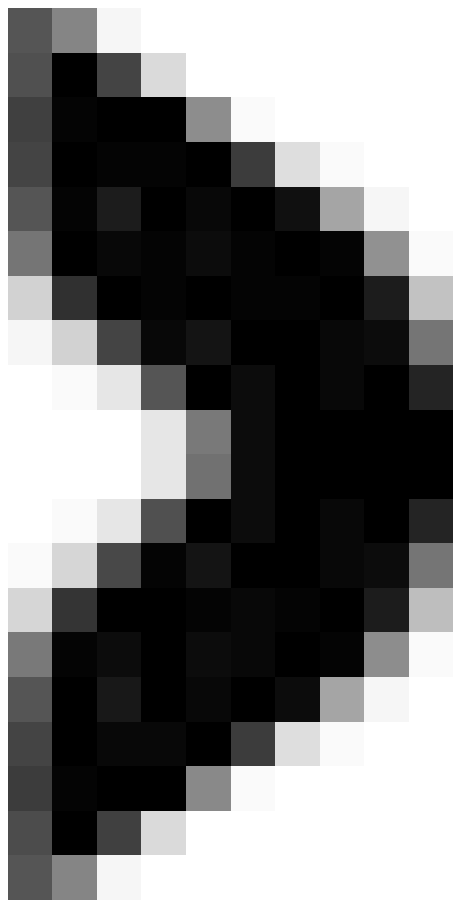}
            \includegraphics[width=0.08\textwidth,clip=true,trim=5cm 1cm 4cm 1cm]{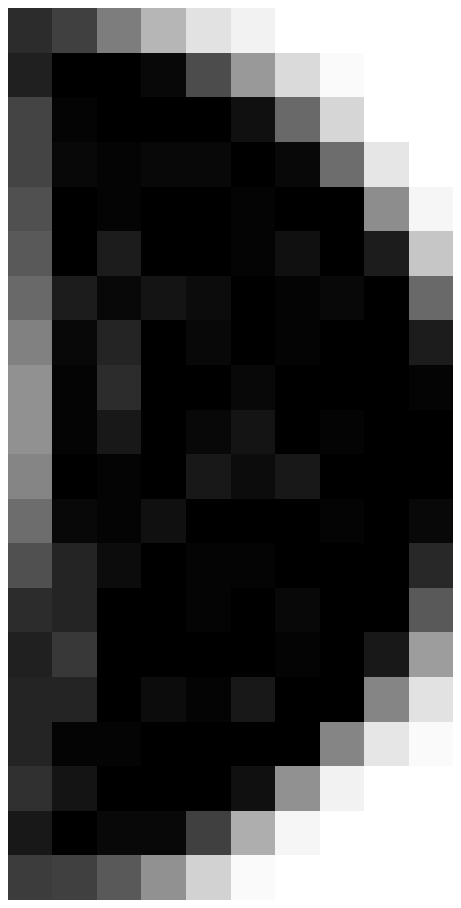}\\
            \includegraphics[width=0.08\textwidth,clip=true,trim=5cm 1cm 4cm 1cm]{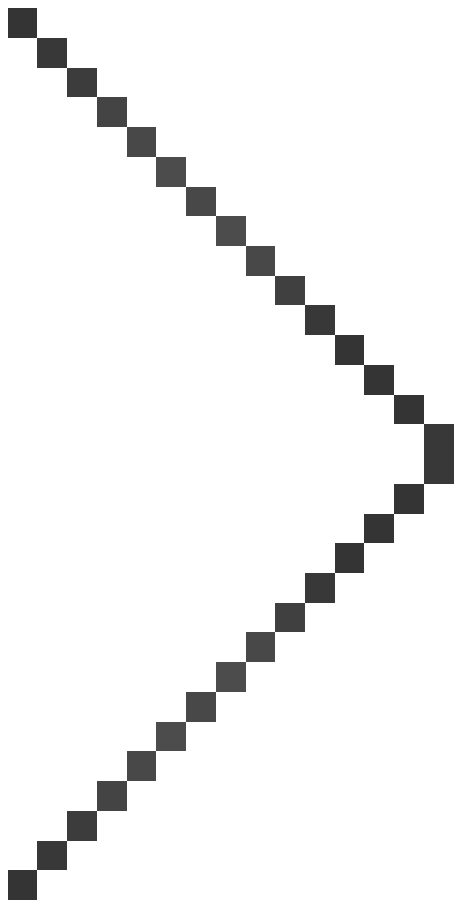}
            \includegraphics[width=0.08\textwidth,clip=true,trim=5cm 1cm 4cm 1cm]{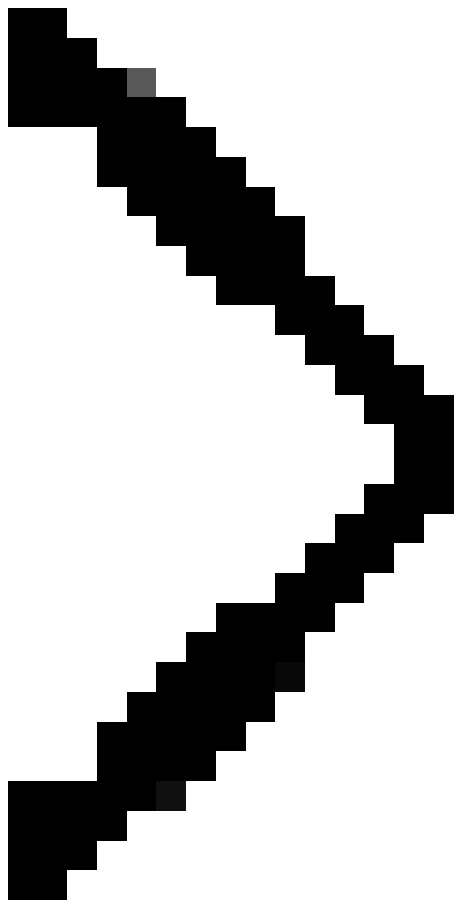}
            \includegraphics[width=0.08\textwidth,clip=true,trim=5cm 1cm 4cm 1cm]{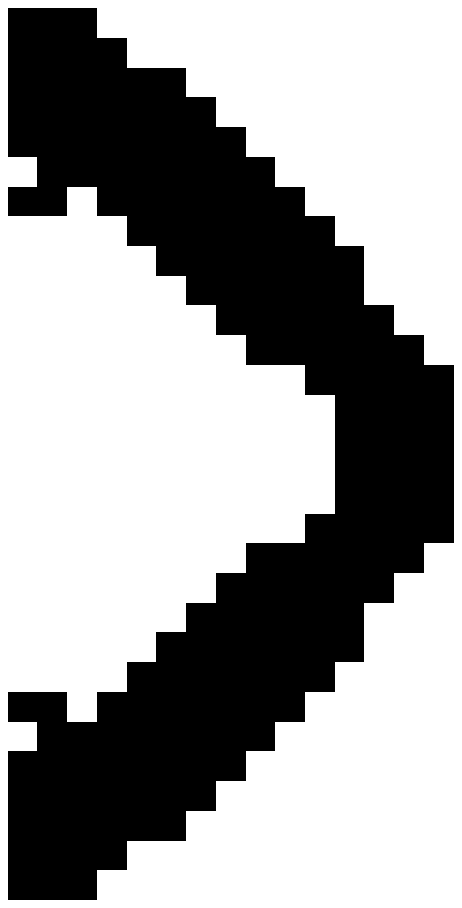}
            \includegraphics[width=0.08\textwidth,clip=true,trim=5cm 1cm 4cm 1cm]{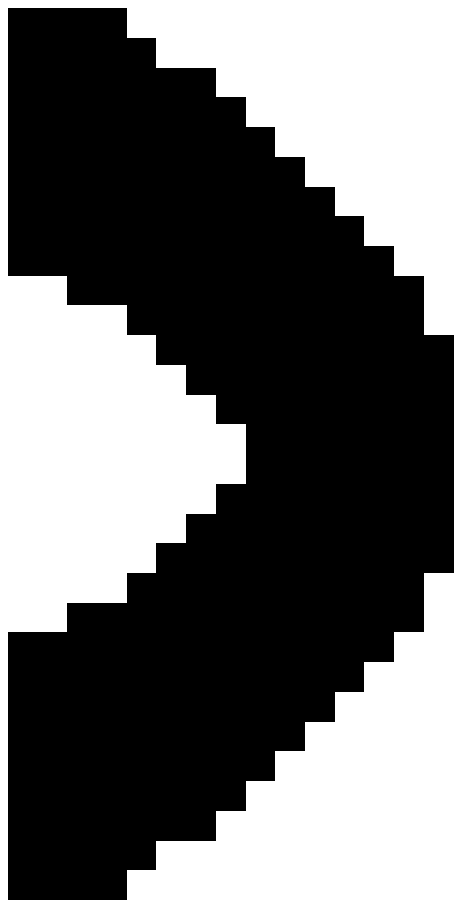}
            \includegraphics[width=0.08\textwidth,clip=true,trim=5cm 1cm 4cm 1cm]{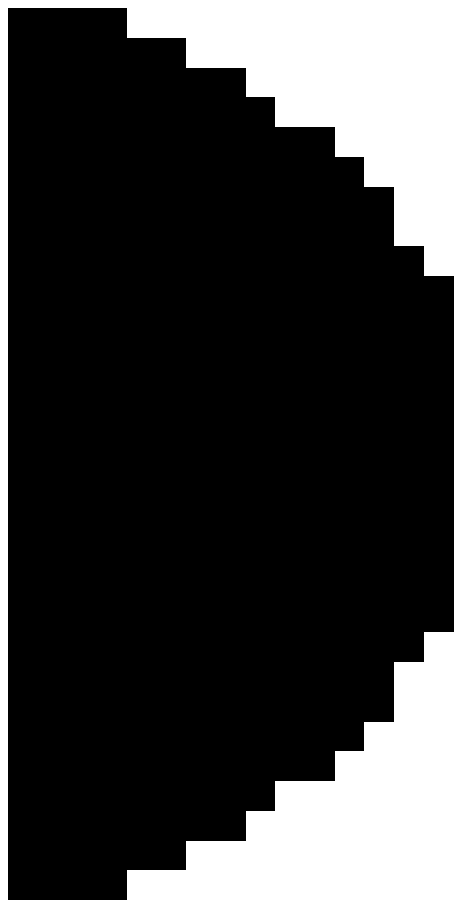}\\
            \vspace{-0.25cm}
            \caption{Density distributions for High Fidelity meshes: EFG-based MNA (top), FEM-based MNA (middle) and SIMP (bottom) for volume fractions (from left to right) 0.05, 0.2, 0.4, 0.6 and 0.8}
            \label{fig:density2}
        \end{figure}

        The figure \ref{fig:vfrac} shows the results: the decreasing efficiency of adding matter as the absolute value of the compliance slope diminishes and the compliance values given by the SIMP smaller that the ones given by the MNA. This can be understood with figure \ref{fig:density2}. It is clear that the material distribution given by the MNA is more blurry and therefore more spread than the one given by the SIMP. The regions where density plays its most important role in global stiffness cannot be filled as much. This is merely due to the small number of structural members used. In fact, smaller compliances can be reached with a finer material description, but not as small as the ones given by the SIMP with the same discretization.
        
        The MNA converges well when the volume fraction is intermediate (from 0.2 to 0.8 approximately). This is also due to the low number of material variables. Low volume fractions cause the structural member to become very thin and therefore not well recognized by the discretization. When the volume fraction is high, the structural members have a tendency to get close in regions where stiffness is needed, not filling completely the designable space. This is why the final volume fraction does not exceed 0.8 in the test case.
        
    \subsection{Filtering}
The effect of a filter is now investigated. Mass nodes will be used as they provide a more intuitive insight of what is happening. The FEM-based MNA is used with a mesh of $10\times20$ elements with $6\times6$ mass nodes.
        \begin{figure}[ht]
            \centering
            \includegraphics[width=0.4\textwidth]{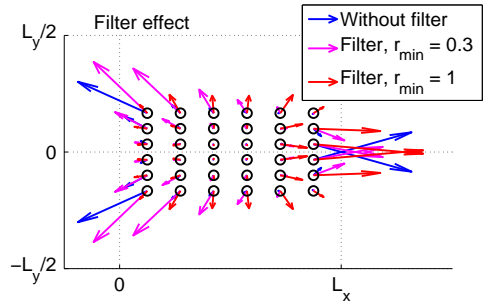}
            \caption{Effect of a filter on the initial compliance sensitivities (the arrows point to decreasing compliance configurations)}
            \label{fig:filterEffect}
        \end{figure}

\subsection{Very first results on L-shape}
\begin{figure}[ht]
    \centering
    \includegraphics[width=0.25\textwidth]{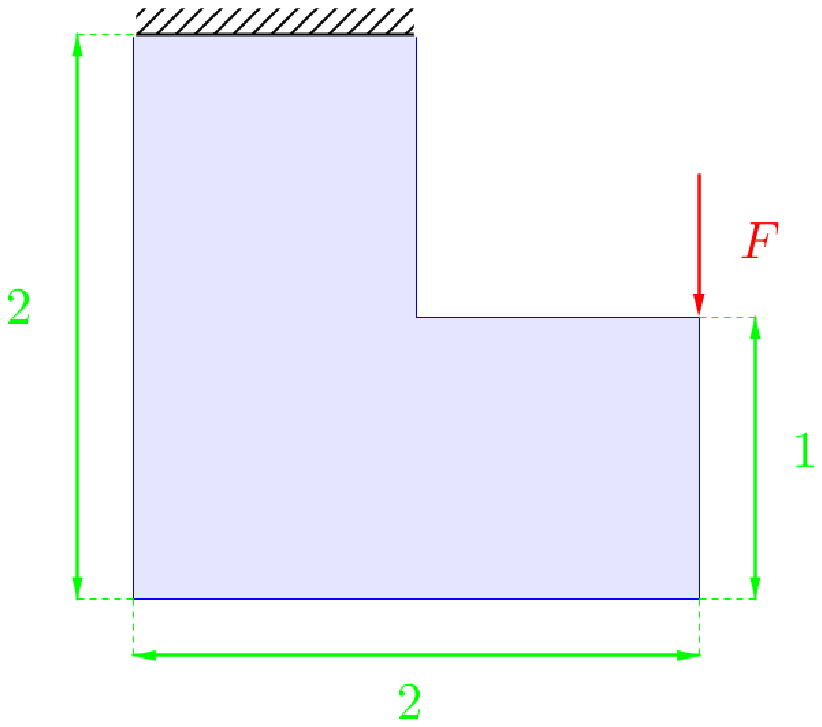}
    \caption{L-shape test case}
    \label{fig:TestCaseLShape}
\end{figure}

Another interesting test case is the L-shape problem described in the figure \ref{fig:TestCaseLShape}. 

Almost all the numerical parameters used are the same as in table \ref{tab:testCase}. If not, they are listed in table \ref{tab:LShape}, which also includes other parameters linked to the merging and suppression of nodes.
    
    \begin{table}[ht]
            \footnotesize
            \def\arraystretch{1.2}
            \centering
            \begin{tabular}{|c|c|c|}
                \hline
                \textbf{Quantity} & \textbf{Symbol} & \textbf{Value} \\
                \hline
                Elements/cells per unit length &$n_x$ or $n_y$ & 10 \\
                \hline
                Num. of mass nodes & $n_{mn}$ & 40 \\
                \hline
                Volume fraction & $v_{Frac}$ & 0.5 \\
                \hline
                Maximum step norm & / & $2/5$ \\
                \hline
                Tolerance on angles & $tol_\theta$ & 5$^\circ$\\
                \hline 
                Tolerance on distances & $tol_d$ & 0.1 \\
                \hline
                Tolerance on dimensions & $tol_l$ & 0.25\\
                \hline
                Denisty radius & $r_\rho$& 0.37\\
                \hline
            \end{tabular}
            \caption{L Shape numerical values}
            \label{tab:LShape}
        \end{table}

Though relatively simple, it raises some issues. First, the Matlab's optimizers \texttt{fminunc} or \texttt{fmicon} appear to be much less efficient in that case. At the first iteration, the nodes with the highest sensitivities are often moved in a awkward configuration, which ends up in making the optimized part not only located at a local minimum of compliance, but also in a configuration which is impossible to manufacture. The problem originates from the optimal step used in the line search. Instead, it seems that using steps with maximum allowable length is more cautious (at the cost of evaluating the sensitivites more often). Hence, a simple steepest descent method has been used in the following. The FEM-based MNA results obtained on a mesh with 10 elements per unit length and $8\times5$ mass nodes are given hereafter.

\begin{figure}[ht]
        \centering
        \includegraphics[width=0.2\textwidth,clip=true,trim=0cm 0cm 1cm 0cm]{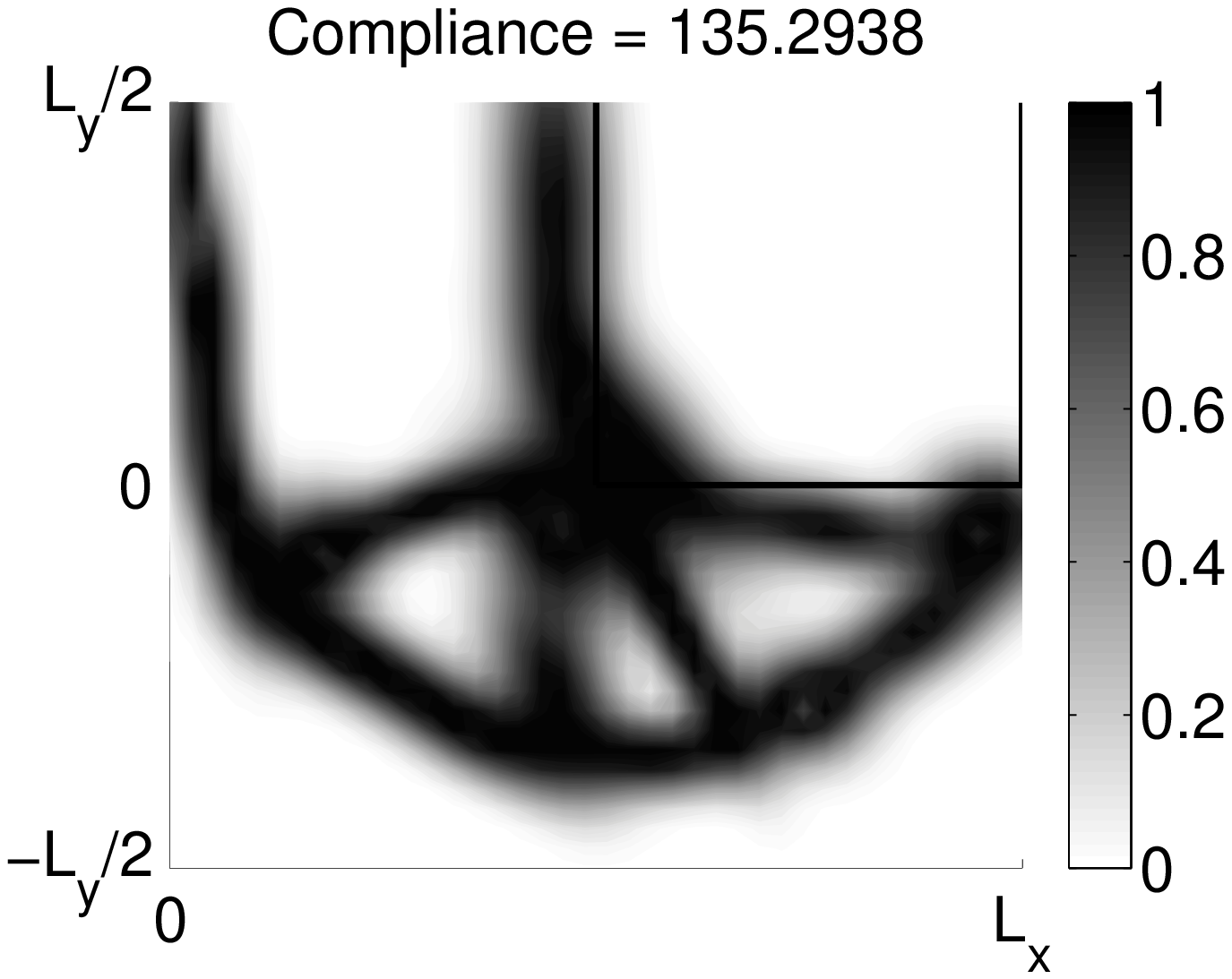}
        \includegraphics[width=0.2\textwidth,clip=true,trim=0cm 0cm 1cm 0cm]{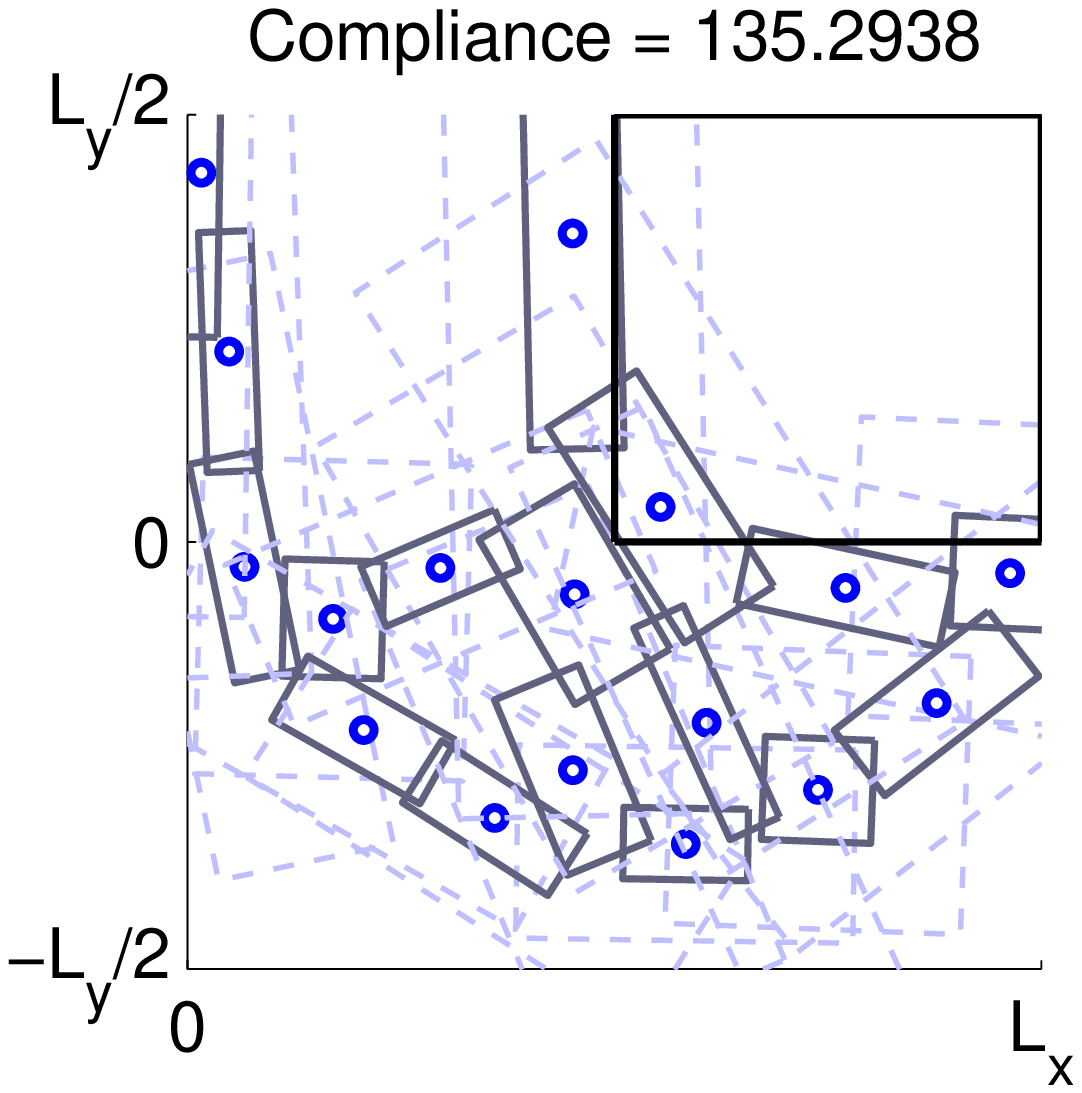}
        \caption{Optimized layout for the L-shape problem with deformable structural members for a volume fractions of 0.5 with the FEM-based MNA with merging and suppression of nodes}
        \label{fig:LShapeMerging}
    \end{figure}

The merging algorithm presented in section \ref{sec:methodology} can then be used. The isolated zero-width nodes, which do not contribute to the part's stiffness, should also be suppressed. Figure \ref{fig:LShapeMerging} presents the results of an optimization taking into account the merging process and the isolated nodes suppression.
     
    The merging and suppression procedures allow the number of nodes to be reduced to 17. The structure appears simpler. It is however a little more compliant. One can note that a structural member lies outside the design domain. This is not really an issue, since the part that lies outside the  domain is not taken into account for the structural stiffness. An hypothetical conversion to a CAD model would however need to eliminate the part of this structural member that lies outside the design domain by a boolean operation for instance.

\section{Toward element recognition}

The interested reader can fin here the 2016 topMNA description (for educational and research purpose). This code is largely inspired from top88.m. Two test cases are illustrated hereafter.

\subsection{Cantilever beam}

\begin{figure}[ht]
    \centering
\includegraphics[totalheight=0.7in]{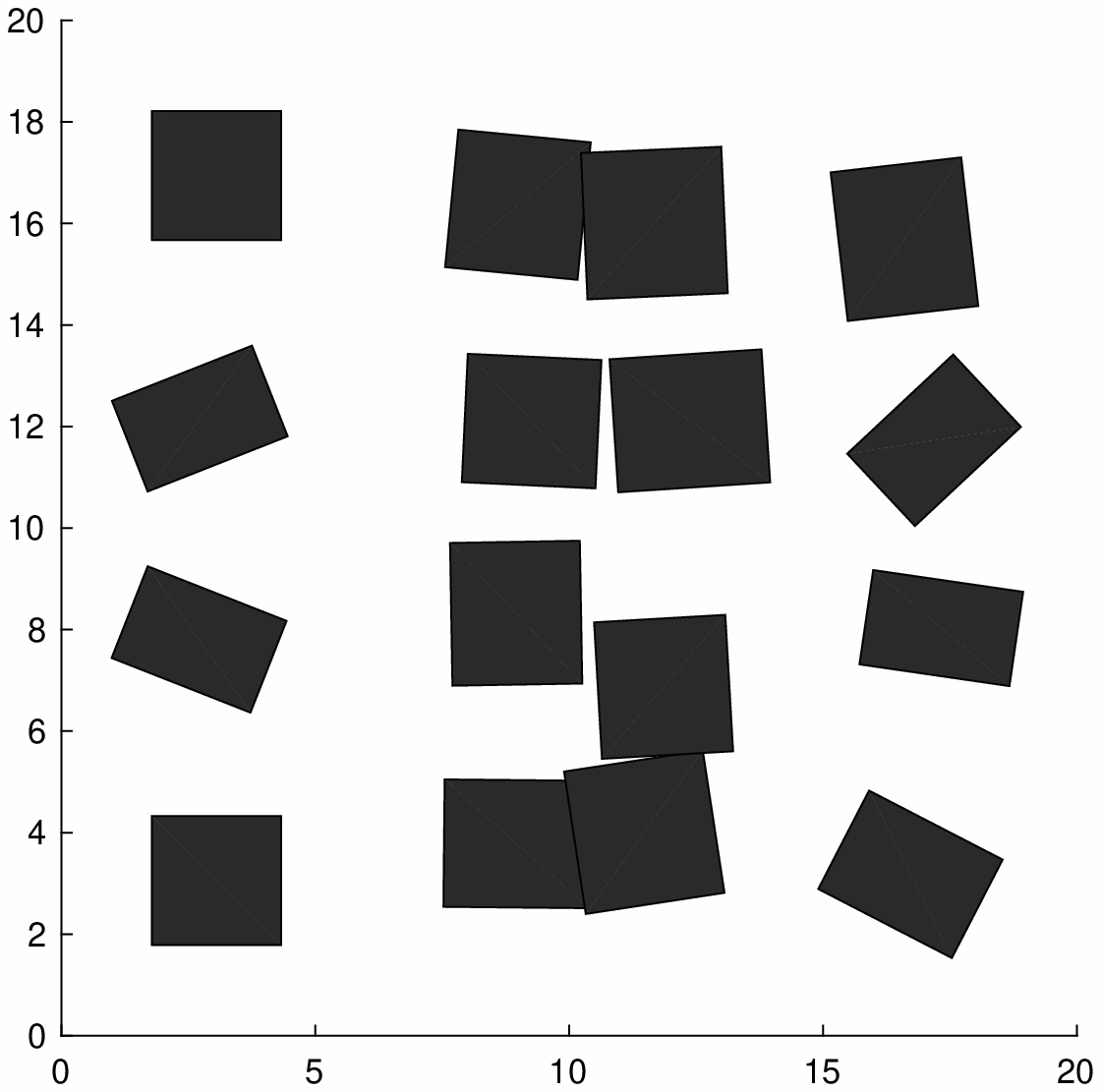} \includegraphics[totalheight=0.7in]{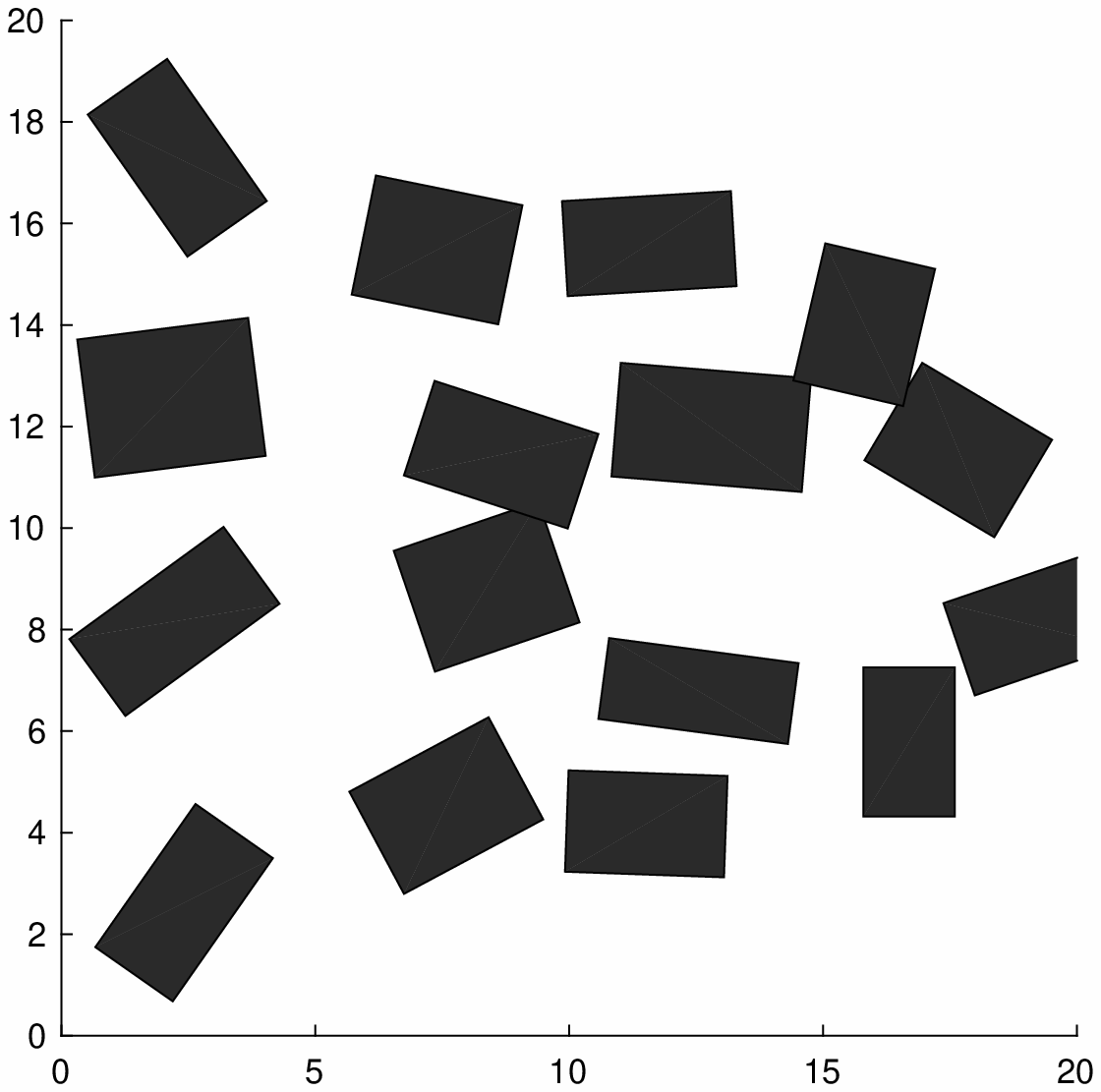}
\includegraphics[totalheight=0.7in]{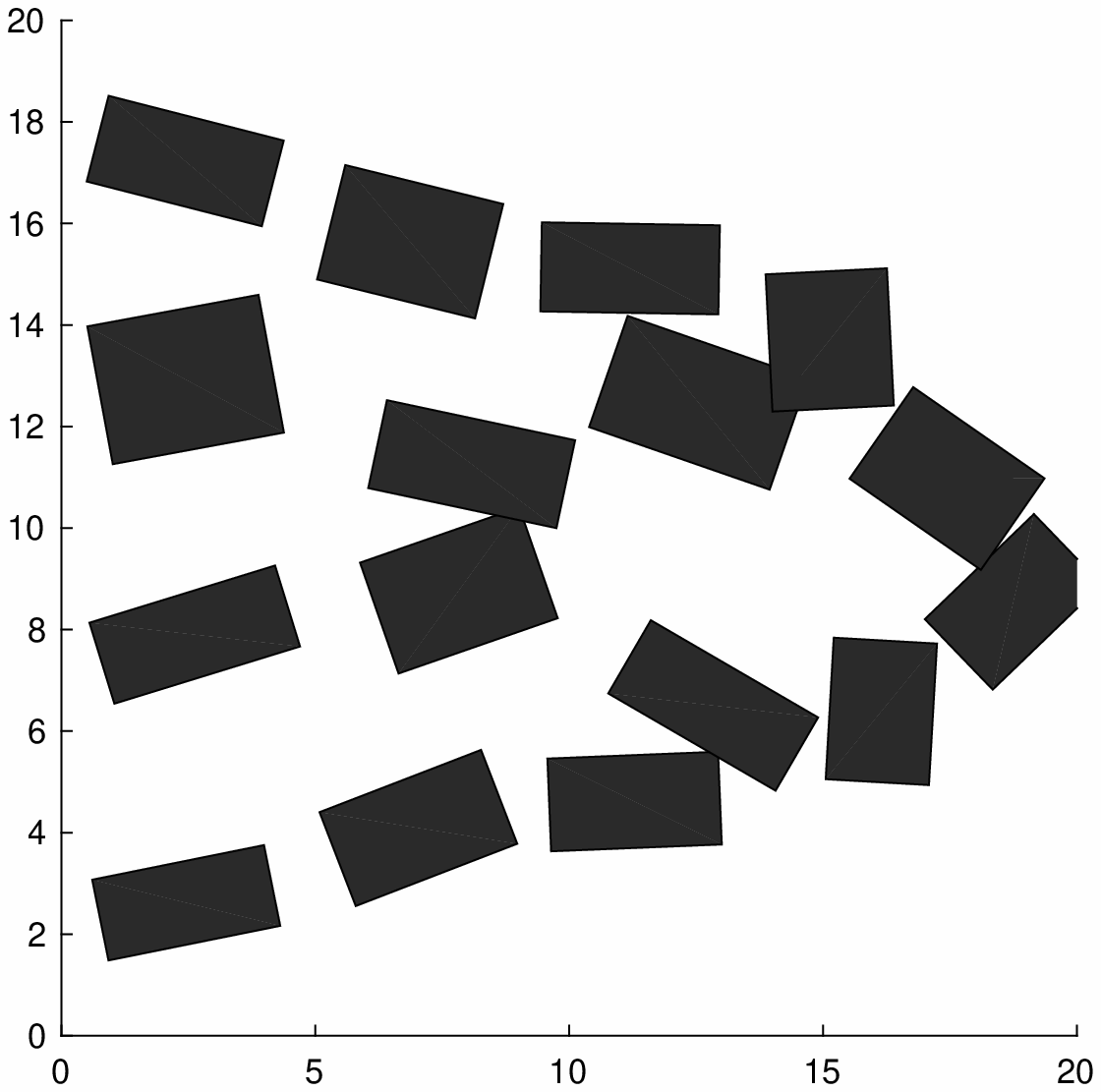}
\includegraphics[totalheight=0.7in]{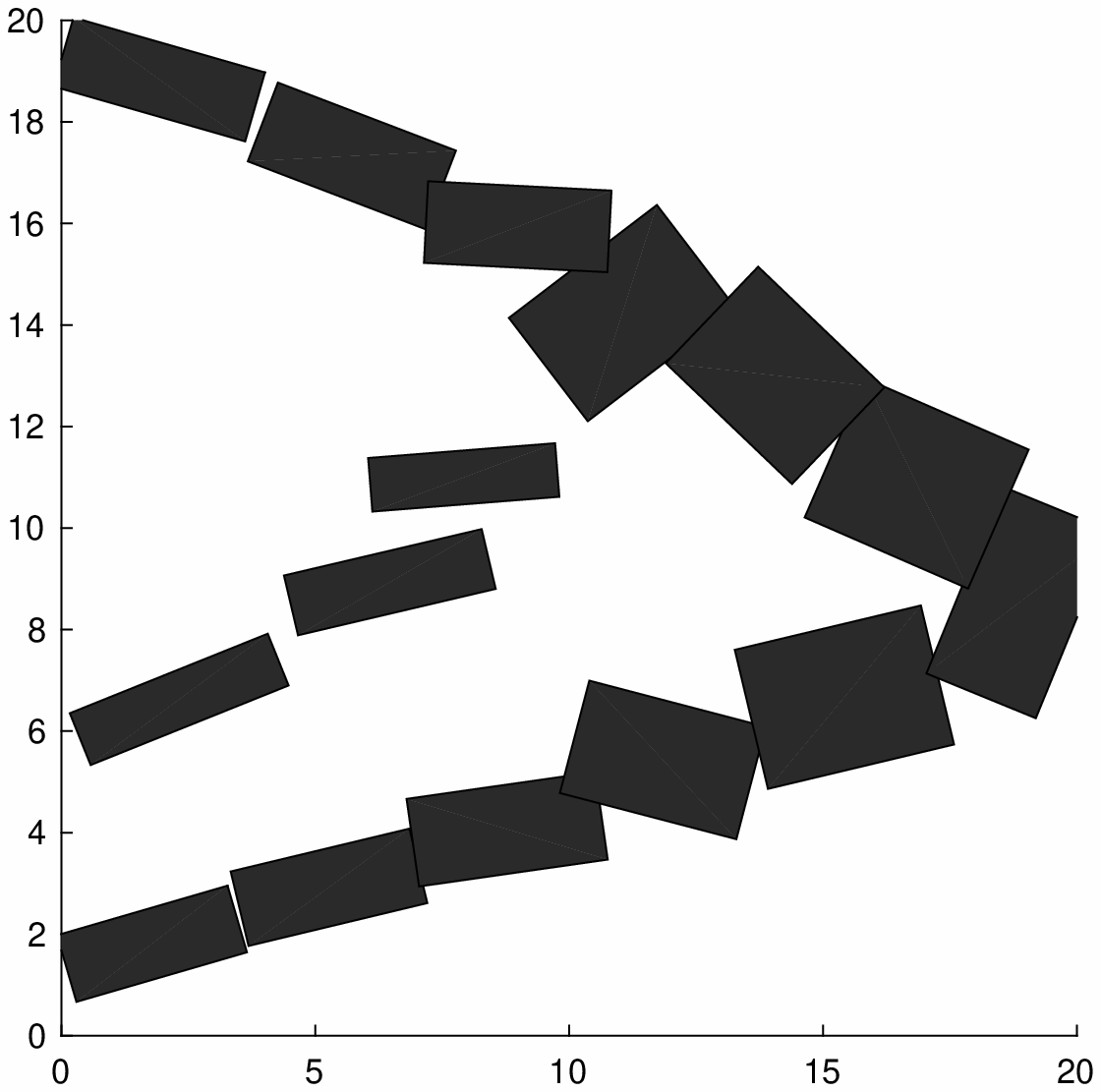}
\includegraphics[totalheight=0.7in]{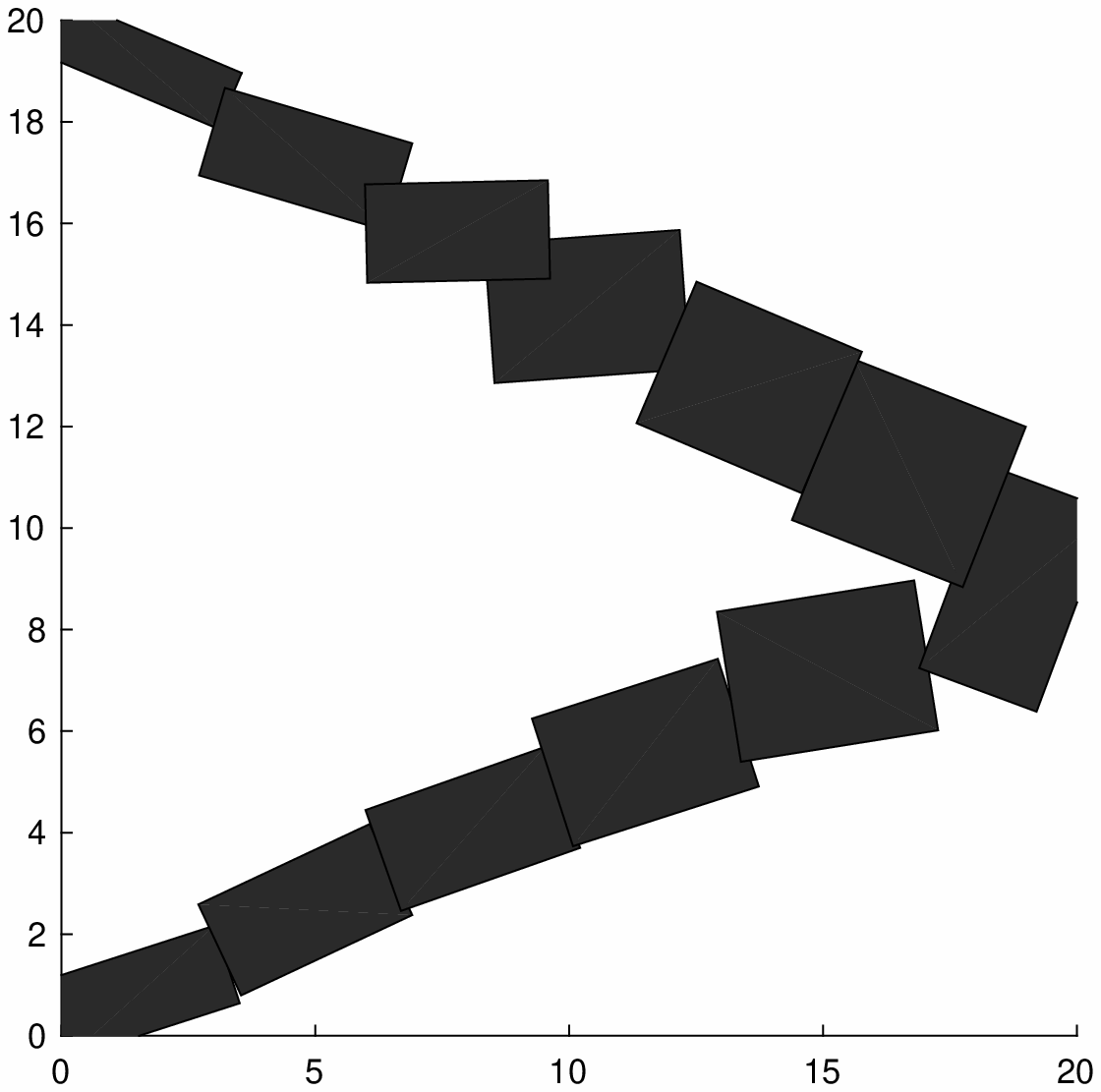}
\includegraphics[totalheight=0.7in]{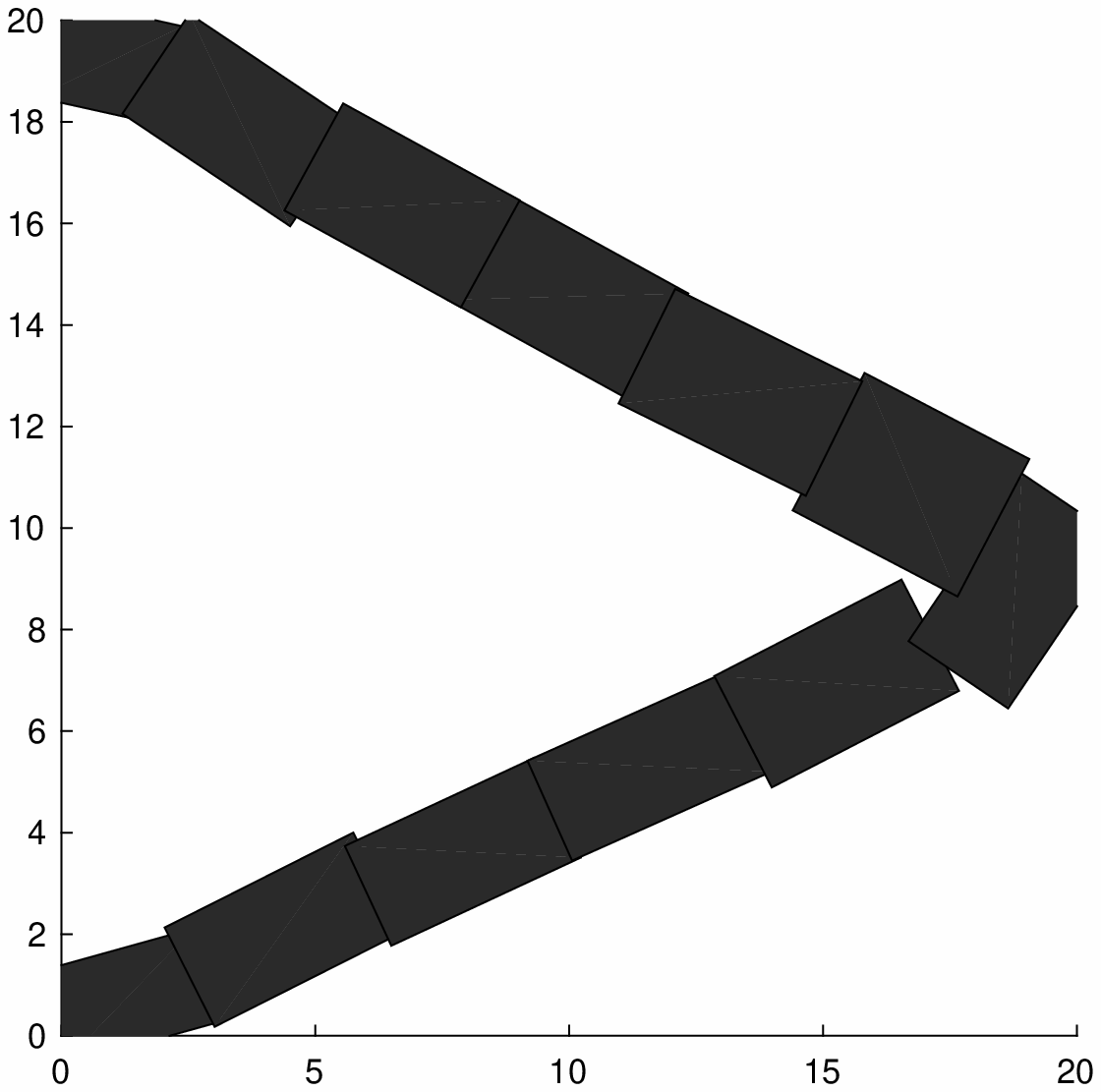}
\includegraphics[totalheight=0.7in]{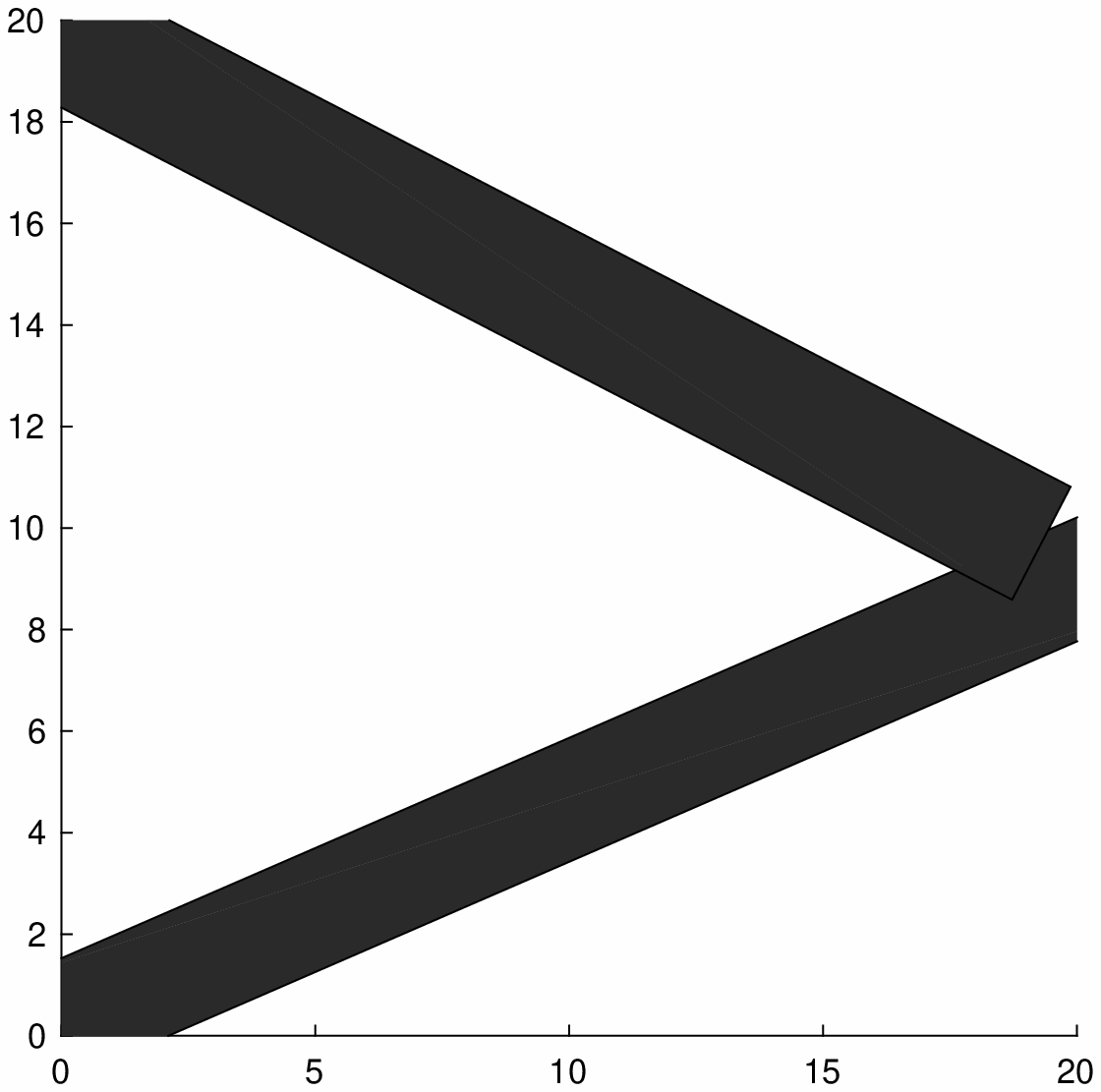}
\includegraphics[totalheight=0.7in]{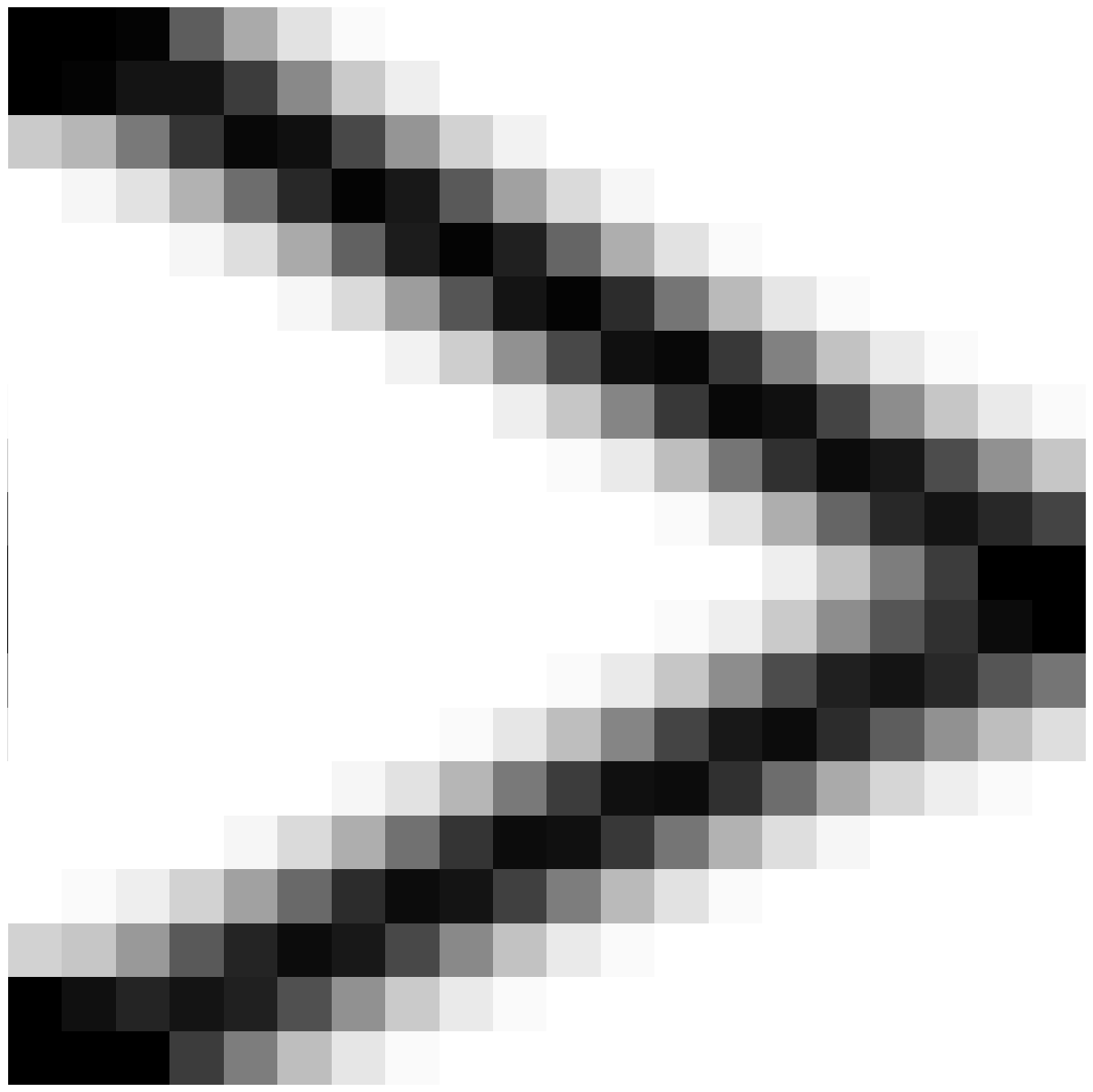}
   \caption{Cantilever beam results: iteration 1,3, 7, 34, 47, 107, 108 on 108 and SIMP}
    \label{fig:Test1}
\end{figure}

We also check the testcase with $volfrac=0.5$ according to \cite{yi2017identifying} with an initial good nodes placement. The optimizer converges locally then oscillates between two solutions.
See file : CheckCantileverMNA.m

\subsection{L-Shape}

\begin{figure}[ht]
    \centering
\includegraphics[totalheight=0.7in]{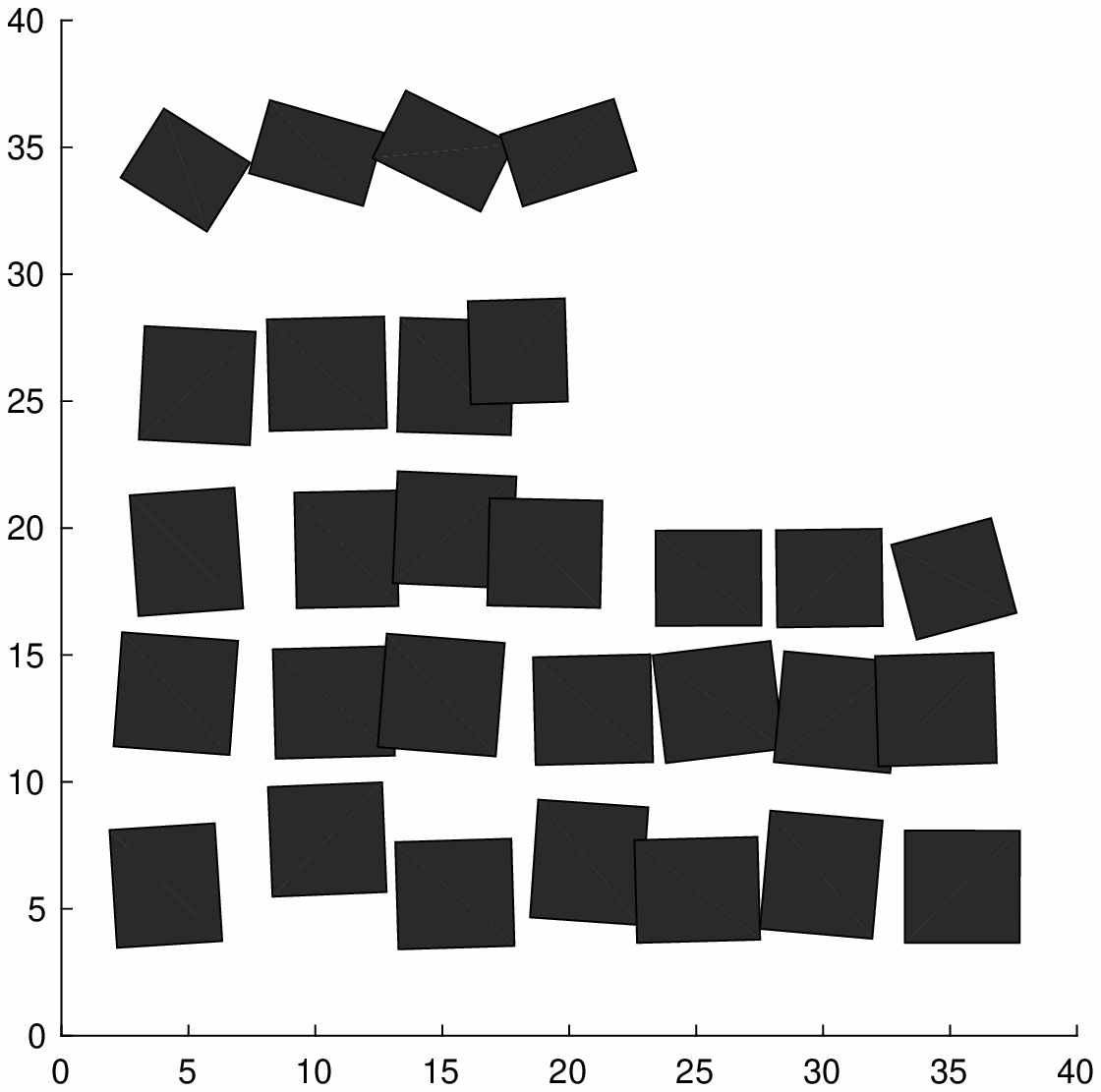} 
\includegraphics[totalheight=0.7in]{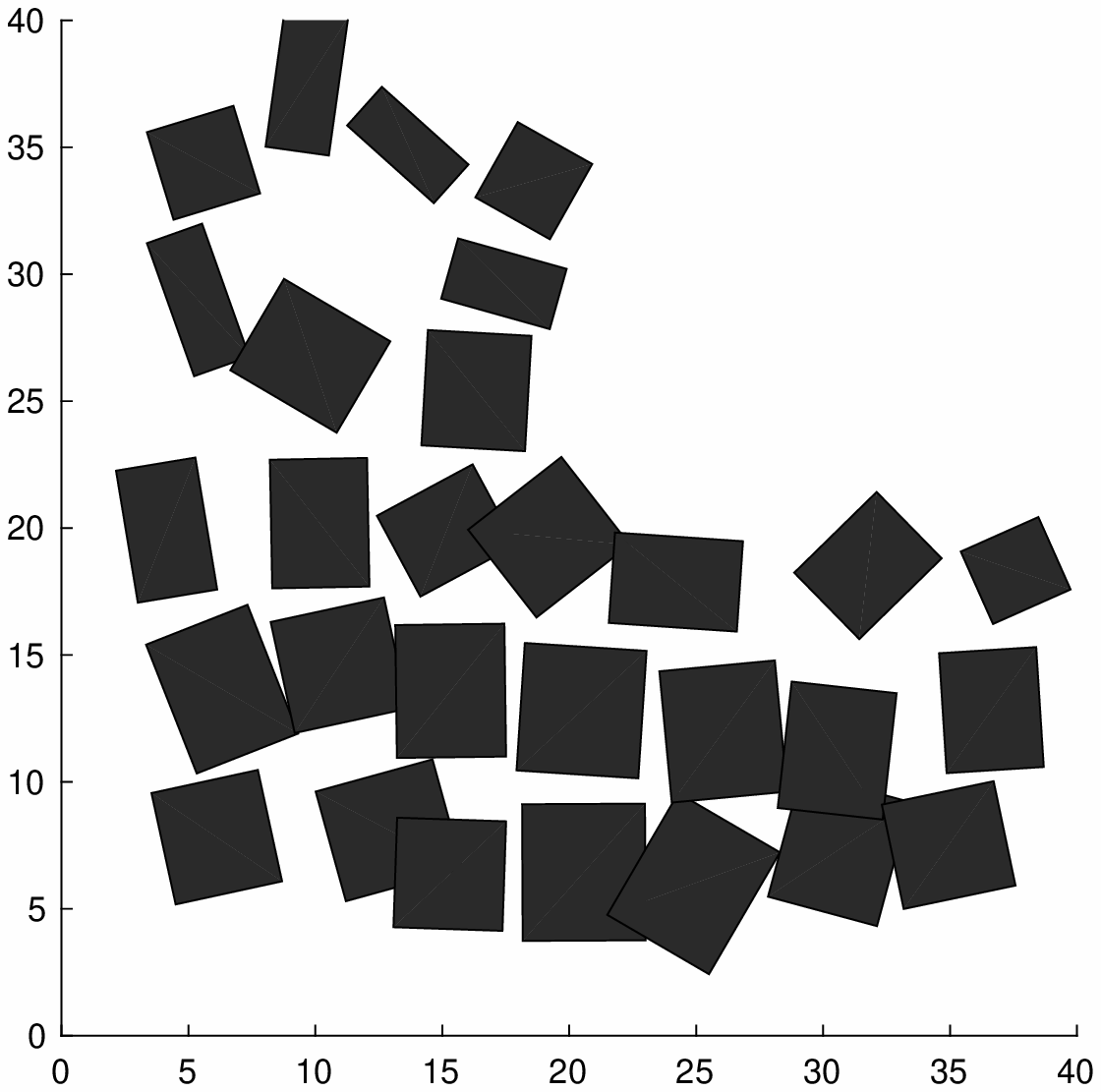}
\includegraphics[totalheight=0.7in]{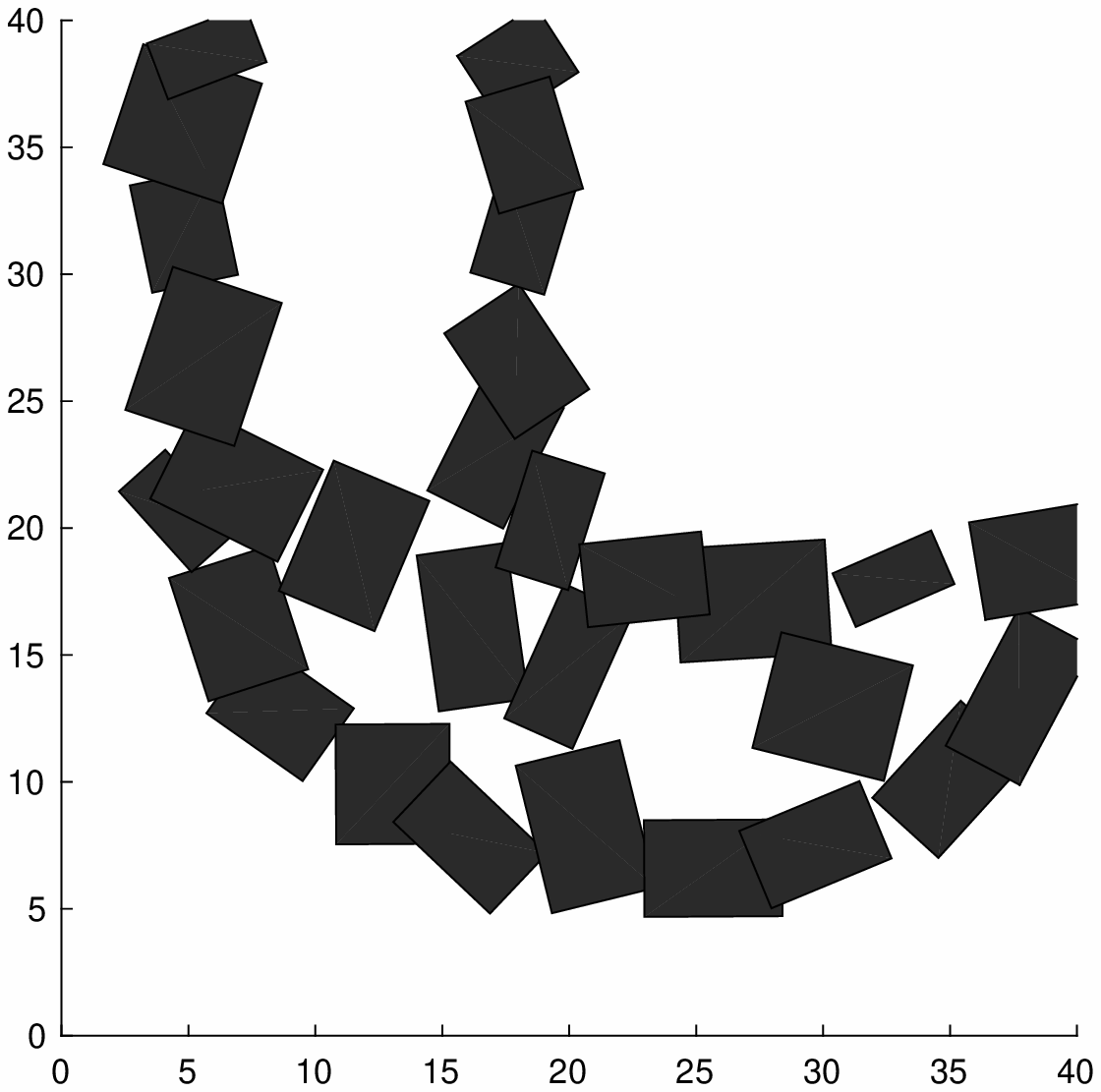}
\includegraphics[totalheight=0.7in]{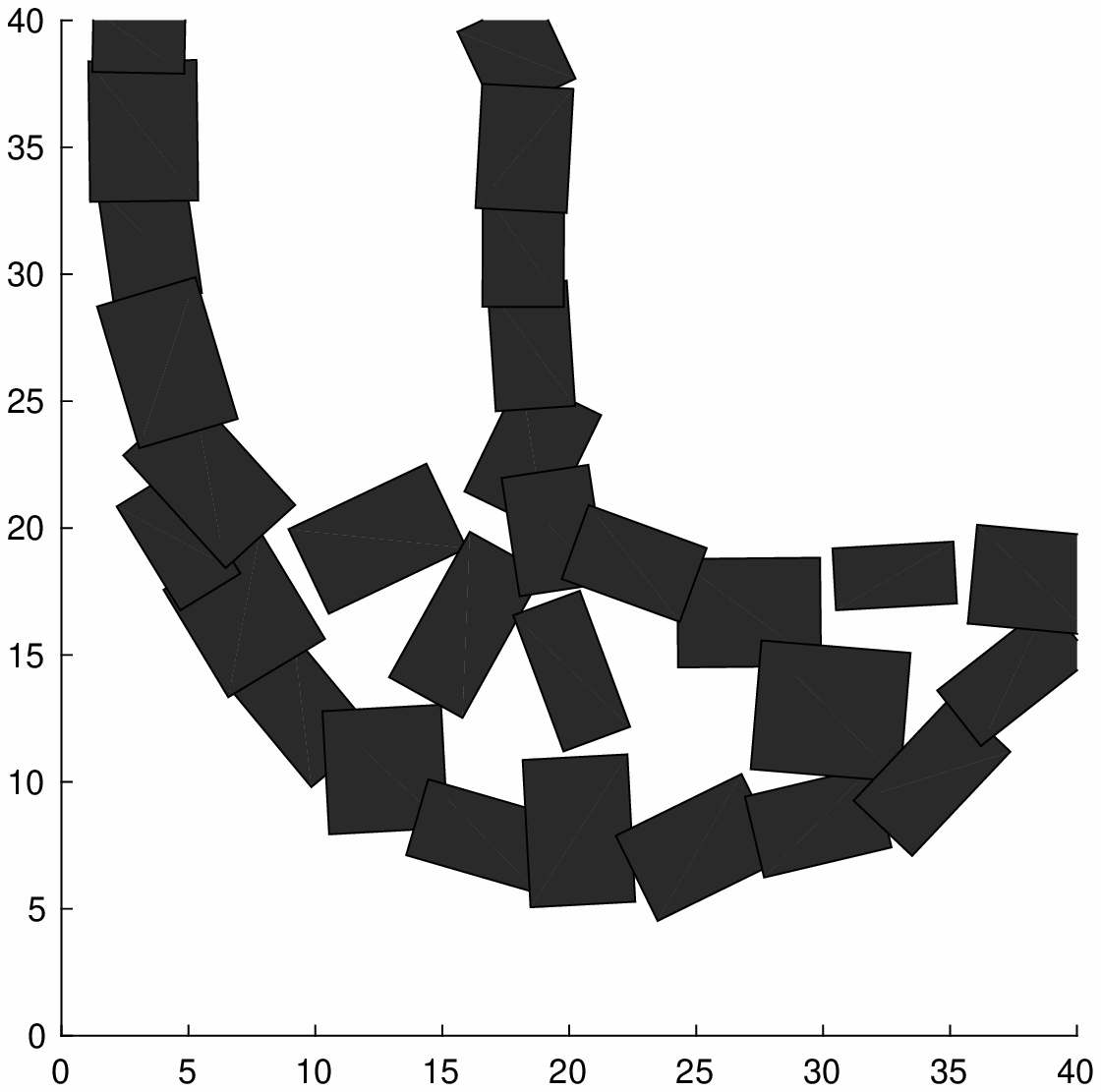}
\includegraphics[totalheight=0.7in]{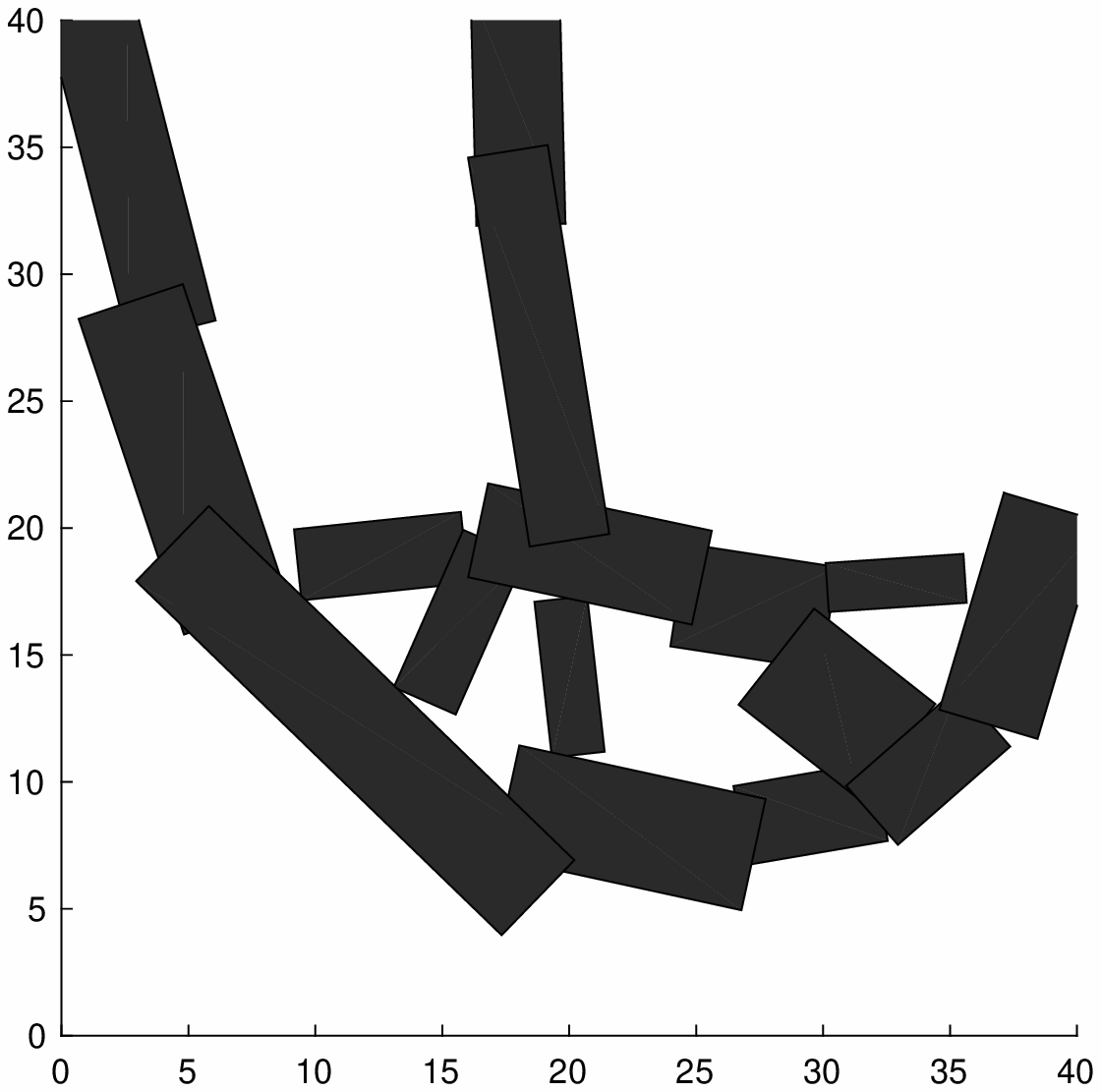}
\includegraphics[totalheight=0.7in]{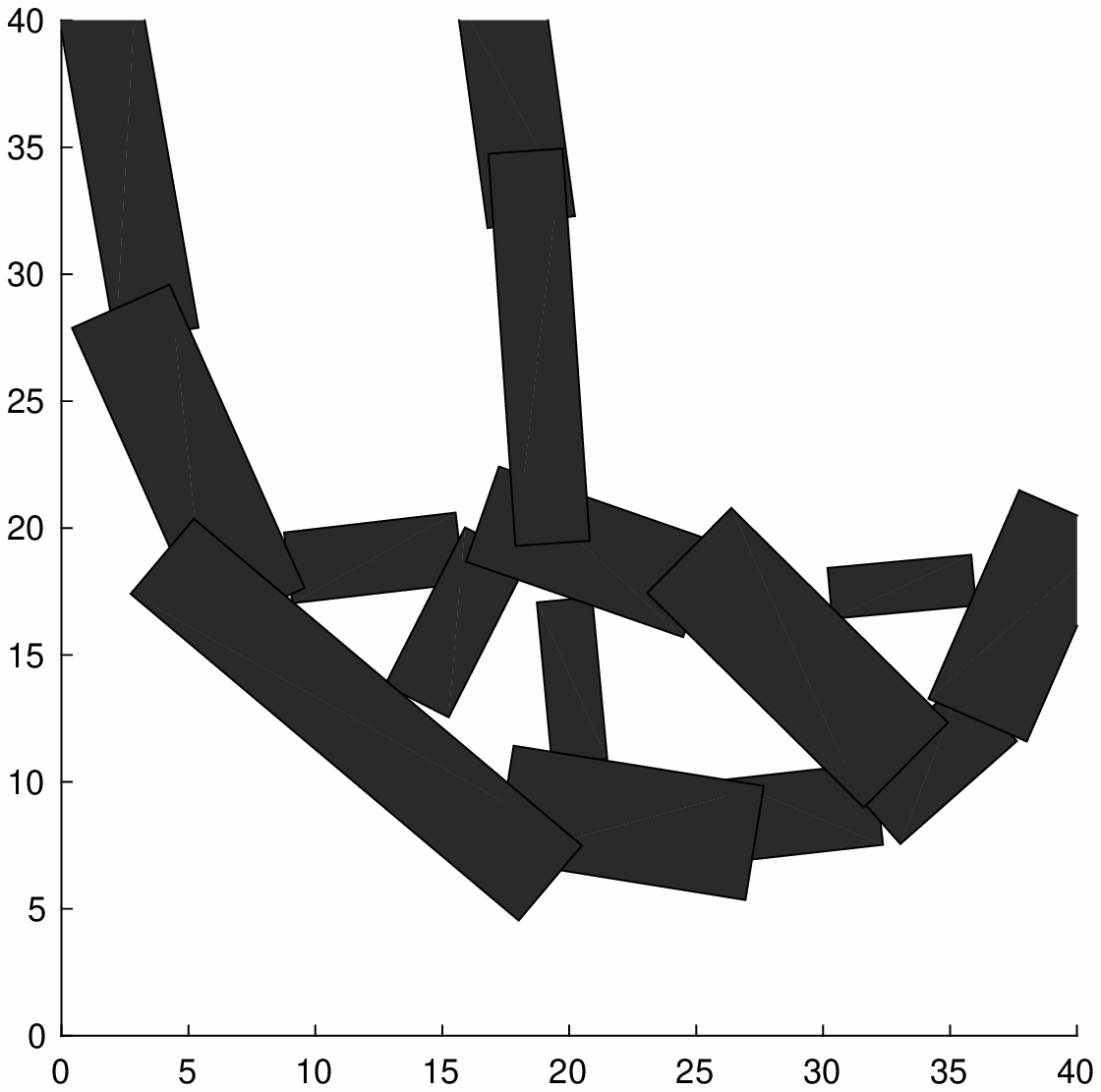}
\includegraphics[totalheight=0.7in]{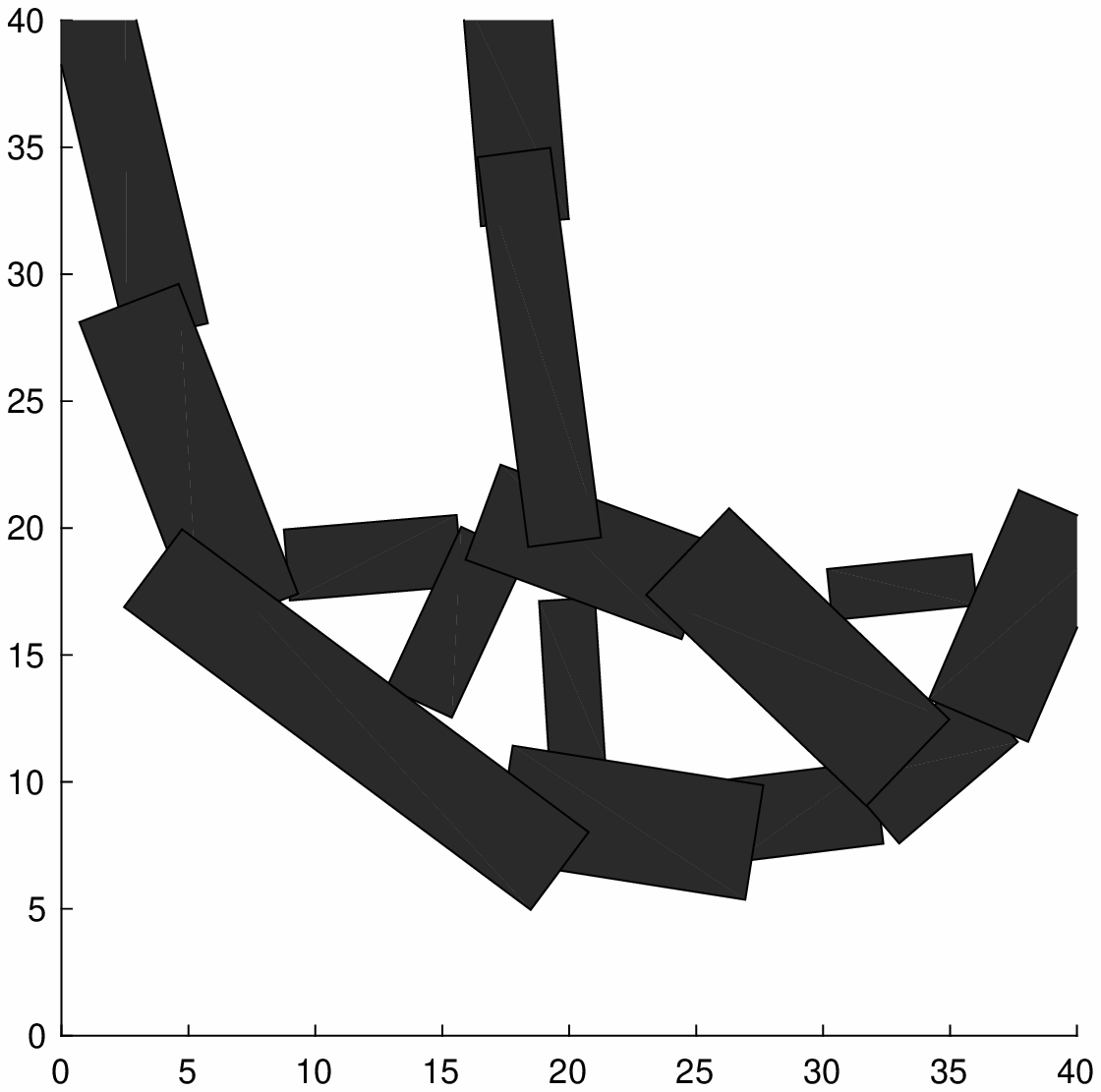}
\includegraphics[totalheight=0.7in]{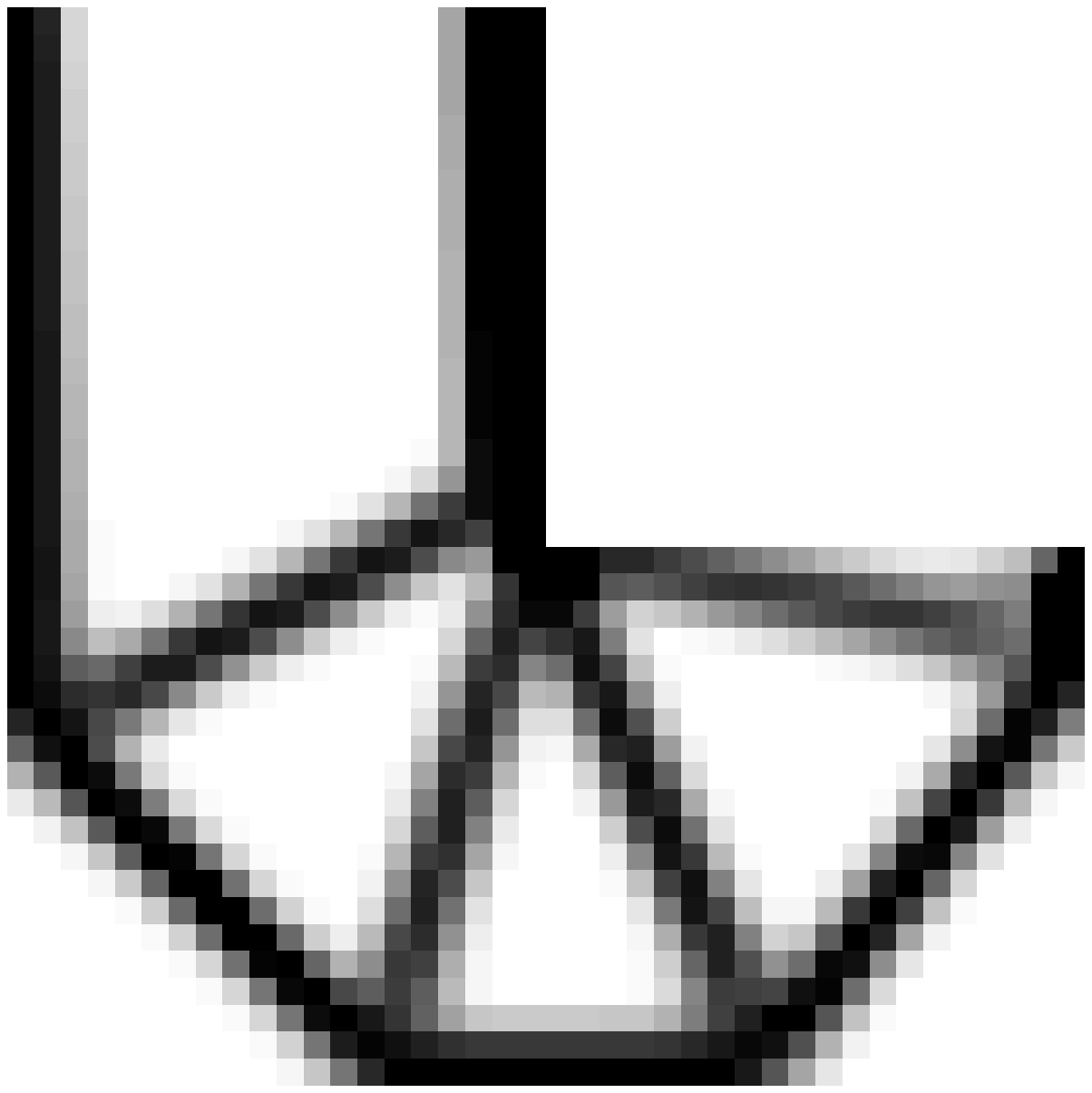}
   \caption{L-shape results: iteration 1,3, 14, 21, 32, 38, 39 on 39 and SIMP}
    \label{fig:Test2}
\end{figure}

\subsection{Ongoing works}
Of course there is 2 hyperparameters (very sensitive) to tune: ratio, aspect.
The first is a tolerance on the distances between the centers of the nodes relatively to their dimensions, the second is a maximum difference of aspect ratio between adjacent elements.
The algorithm loops on a tolerance on the change on design variable (change).

We still have several items to treat:
\begin{itemize}
    \item Our gradient based method is of course dependent of x0. But x size is changing during optimization due to elimination of disconnected nodes
    \item Our optimizer is not the best: we for sure reach local minimum (see L-shape results). We currently use MMA (or derived algorithms) much more adapted to structural optimization \cite{svanberg1987method}
    \item Our method is stabilized when the penalization is increasing step by step
    \item Our merging criteria is definitely not perfect
\end{itemize}

Our full developments (since 2016) can be found in:
\url{https://github.com/GhislainRaze/Topology-Optimization}

The proposed following examples in: 
\url{https://github.com/jomorlier/Topology-Optimization/tree/master/arXiv}

\section{Conclusions}\label{Section:Conclusions}

 The MNA is an interesting approach for topology optimization. Compared to the SIMP, it is computationally more expensive and requires more effort to be implemented. However, it allows control over the final design by setting the initial number of structural elements and the bounds of their characteristics. During the optimization, the structural elements can be merged or suppressed according to tolerances provided by the user. It also yields designs which are far easier to interpret: positions, orientations and dimensions of structural members. Namely, the conversion of the output data to CAD models for the shape optimization or ALM models could be much more direct than with the other topology optimization techniques with an efficient shape recognition algorithm, progressing towards design automation.
    
Using deformable structural members is particularly advantageous with this approach as it greatly reduces the number of design variables. Their change in dimensions can be troublesome for the discretization but solutions can be found. For instance, minimum dimensions for the structural members and constraints of the optimizer can be set. A deformable structural element that would become too small could also be suppressed, leaving the available material to other more important structural elements.
    
The adaptability of this approach and the associated code to 3D problems is quite straightforward. Computational time could probably be gained by improving the code. Notably, the most time-consuming part of the MNA is the evaluation of the density field. Since it is done independently on each Gauss point, it could be interesting to use parallel computing (and GPU) to gain time.

\section*{Acknowledgement}

\bibliographystyle{plainnat}
\bibliography{bibliography.bib}

\end{document}